\documentclass[11pt]{article}
\usepackage[utf8]{inputenc}

\usepackage{textcomp}
\usepackage[russian, english]{babel}

\usepackage{amsmath, amsfonts, amssymb,amsthm}
\usepackage{verbatim}

\usepackage{float}
\usepackage{changepage}
\usepackage{array}
\usepackage{enumerate}
\usepackage{hyperref}
\usepackage{authblk}

\usepackage[nottoc,numbib]{tocbibind}

\usepackage{xcolor}

\hypersetup{colorlinks,linkcolor={blue},citecolor={red},urlcolor={blue}}

\usepackage{tikz}
\usetikzlibrary{calc, matrix, arrows,decorations.pathmorphing, positioning}
\numberwithin{equation}{section}
\theoremstyle{definition}

\newtheorem{conjecture}{Conjecture}[section]
\newtheorem{theorem}{Theorem}[section]

\newtheorem{definition}{Definition}[section]
\newtheorem{remark}{Remark}[section]

\usepackage{etoolbox}
\makeatletter
\AfterEndEnvironment{definition}{\@doendpe}
\AfterEndEnvironment{theorem}{\@doendpe}
\AfterEndEnvironment{lemma}{\@doendpe}
\AfterEndEnvironment{prop}{\@doendpe}
\AfterEndEnvironment{remark}{\@doendpe}
\AfterEndEnvironment{example}{\@doendpe}
\AfterEndEnvironment{consequence}{\@doendpe}
\makeatother

\newcommand{\beq}{\begin{equation}}
\newcommand{\ee}{\end{equation}}

\newcommand{\ben}{\begin{equation*}}
\newcommand{\een}{\end{equation*}}

\usepackage{caption}

\DeclareMathOperator{\Fib}{Fib}

\DeclareMathOperator*{\Vac}{Vac}

\DeclareMathOperator*{\Tr}{Tr}

\usepackage{paralist}

  \pltopsep=\medskipamount
\newcommand{\qbinom}[2]{\genfrac{[}{]}{0pt}{}{#1}{#2}}

\usepackage[toc]{appendix}

\providecommand{\keywords}[1]
{
  \small	
  \textbf{\textit{Keywords---}} #1
}

\newcolumntype{M}[1]{>{\centering\arraybackslash}m{#1}}

\title{Durfee rectangle identities as character identities for infinite Fibonacci configurations}

\author[1]{Timur Kenzhaev \thanks{kenzhaev\_t\_d@mail.ru}}

\affil[1]{Skolkovo Institute of Science and Technology, Moscow, Russia}

\date{}

\begin{document}

\pagenumbering{arabic}

\maketitle

\begin{abstract}
We introduce a natural generalization of Maya diagrams~---~the space of infinite Fibonacci configurations, which are specified functions on $\mathbb{Z}$ with values $1$ and $0$. Infinite Fibonacci configurations are particularly interesting as soon as they parametrize Feigin-Stoyanovsky type bases in lattice vertex superalgebras $V_{\sqrt{N}\mathbb{Z}}$ and their irreducible modules. We calculate the character of such configurations space by two different ways and obtain series of combinatorial identities. These identities turn out to be Durfee rectangle identities with shifts in base and height. 
\end{abstract}

\begin{center}
\section{Introduction}
\end{center}
This article is devoted to the character identities coming from Feigin-Stoyanovsky bases of one-dimensional lattice VOAs. In the pioneering work \cite{FS} Feigin and Stoyanovsky built the semi-infinite basis of $\widehat{\mathfrak{sl}_2}$ standard modules $L_{(0, 1)}$ and $L_{(1, 1)}$. Namely, let $e, f$ and $h$ be a standard basis in $\mathfrak{sl}_2$, $e(z) = \sum\limits_{n\in\mathbb{Z}}\,e_n\,z^{- n - 1}$, then the following theorem holds:
\begin{theorem}
\label{TheoremSL2}
$L_{(0 , 1)}$ has the basis of semi-infinite monomials $e_{i_1}\,e_{i_2}\,e_{i_3}\dots$
such that 
\begin{enumerate}
    \item $i_1 < i_2 < i_3 < \cdots$ ;
    \item $i_{k + 1} - i_k \geq 2$;
    \item $i_{k} = 1 \mod 2 \quad \text{for } k \gg 1$;
    \item $i_{k + 1} - i_k = 2 \quad \text{for } k \gg 1$.
\end{enumerate}
$L_{(1, 1)}$ has the basis of semi-infinite monomials $e_{i_1}\,e_{i_2}\,e_{i_3}\dots$
such that
\begin{enumerate}
    \item $i_1 < i_2 < i_3 < \cdots$,
    \item $i_{k + 1} - i_k \geq 2$,
    \item $i_{k} = 0 \mod 2 \quad \text{for } k \gg 1$;
    \item $i_{k + 1} - i_k = 2 \quad \text{for } k \gg 1$.
\end{enumerate}
\end{theorem}
The proof is based on the consideration of so-called \textit{principal subspaces}, which became the subject of intensive study \cite{Baranovi2011, Calinescu2007, CALINESCU2008, Calinescuu2008, Calinescu2010, Trupevi2010, zbMATH05919368}. 
\\\\
Currents $e(z)$ and $f(z)$ can be realised as a bosonic vertex operator $V_{\pm\sqrt{2}}(z)$ corresponding to the lattice VOA $V_{\sqrt{2}\mathbb{Z}}$. Thus, \hyperref[TheoremSL2]{Theorem} \ref{TheoremSL2} states the semi-infinite bases of $V_{\sqrt{2}\mathbb{Z}}$ and its irreducible modules. A natural generalization of this result for an arbitrary integral one-dimensional positive definite lattice was proved in \cite{FF_Principle} and \cite{Penn_Lattice}. Let $V_{\sqrt{N}\mathbb{Z}} = \bigoplus\limits_{m\in\mathbb{Z}} F_{m\sqrt{N}}$ be a vertex superalgebra associated with lattice $\mathbb{Z}\lambda$, $\langle\lambda|\lambda\rangle = N$, take conformal vector $\omega = \frac{a_{-1}^2}{2} + \frac{\sqrt{N}}{2}\,a_{-2}$. $V_{\sqrt{N}\mathbb{Z}}$ is generated by two fields $V_{\pm\lambda}(z)$ with OPE
\ben
V_{\lambda}(z)\,V_{-\lambda}(\omega) = (z - \omega)^{-\langle\lambda|\lambda\rangle}\,:V_{\lambda}(z)V_{-\lambda}(\omega):.	
\een
The set of irreducible $V_{\sqrt{N}\mathbb{Z}}$-modules is
\ben
V_{(i), \sqrt{N}} = \bigoplus\limits_{m\in\mathbb{Z}} F_{\frac{i}{\sqrt{N}} + m\sqrt{N}} \hspace{5mm} i = 0, 1, \ldots, N - 1,
\een
where $F_{\frac{i}{\sqrt{N}} + m\sqrt{N}}$ is Fock module with the highest weight $(1, \frac{i}{\sqrt{N}} + m\sqrt{N})$. 
Then the generalization may be summarized as   
\begin{theorem}
\label{TheoremBasis}
Let $V_{\lambda}(z) = \sum\limits_{n\in\mathbb{Z}} e_n\,z^{-n}$. Then $V_{(i), \sqrt{N}}$ has the basis of semi-infinite monomials $e_{i_1}\,e_{i_2}\,e_{i_3}\dots$
such that
\begin{enumerate}
    \item $i_1 < i_2 < i_3 < \cdots$;
    \item $i_{k + 1} - i_{k} \geq N$;
    \item $i_{k} = i \mod N \quad \text{for } k \gg 1$;
    \item $i_{k + 1} - i_{k} = N \quad \text{for } k \gg 1$.
\end{enumerate}
\end{theorem}  
In the case $N = 1$ (Fermionic Fock space) there is a character identity
\beq
\label{Jacobi_Int}
\sum\limits_{n\in\mathbb{Z}} \frac{z^n\,q^{\frac{n(n - 1)}{2}}}{(q)_{\infty}} =  \prod\limits_{j = 0}^{\infty} (1 + zq^j) \prod\limits_{j = 1}^{\infty} (1 + z^{-1}q^j)
\ee
due to this basis, where
$$
(q)_n = \prod\limits_{i = 1}^{n} (1 - q^i) \hspace{5mm} n\in\mathbb{N}_0 \cup \{\infty\}.
$$
 It is the famous \textit{Jacobi's triple product identity} \cite{Jacobi_Fundamenta}, the left side is the character of the direct sum of Fock modules, the right side is obtained by ``cutting" semi-infinite monomials to the left and right parts. We parametrize bases from \hyperref[TheoremBasis]{Theorem} \ref{TheoremBasis} by natural generalization of Maya diagrams~---~infinite Fibonacci configurations and obtain generalization of Jacobi's triple product identity for every $V_{(i), \sqrt{N}}$. This generalization turns out to be equivalent to the series of Durfee rectangle identities with the natural ratio of height to base $1: N$ with shifts. General family of these identities can be written as
\beq
\label{Final identity l + 1 Int}
\frac{1}{(q)_{\infty}} = \sum\limits_{j = 0}^{m - 1}\frac{q^{nj}}{(q)_{n - 1}(q)_j} + \sum\limits_{k = 0}^{\infty} \frac{q^{(k + n)((l + 1)k + m)}}{(q)_{k + n}(q)_{(l + 1)k + m}} + \sum\limits_{i = 1}^{l }\sum\limits_{k = 0}^{\infty} \frac{q^{(k + n + 1)((l + 1)k + m + i)}}{(q)_{k + n}(q)_{(l + 1)k + m + i}}
\ee
for $l, m, n \in\mathbb{N}_0$ with convention $\frac{1}{(q)_n} = 0$ for $n < 0$. These identities for $n = m = 0$ were studied by G. Andrews in \cite{Andrews_Part1}. Note that \eqref{Jacobi_Int} is equivalent to the series of Durfee rectangle identities (by equating coefficients on degrees of $z$-variable)
\ben
z^s: \hspace{5mm}\frac{1}{(q)_{\infty}} = \sum\limits_{\substack{m, l \geq 0 \\ m - l = s \in \mathbb{N}_0}} \frac{q^{ml}}{(q)_m (q)_l}
\een
with fixed $s\in\mathbb{Z}$. These identities are particular cases of \eqref{Final identity l + 1 Int} with the ratio of height to base equal to $1 : 1$.	
\\\\
The article has the following structure~---~in \hyperref[Section 1]{Section 1} we introduce the space of infinite Fibonacci configurations $\Fib_{\infty}^{(\theta,\, l)}$ and calculate its character. \hyperref[Section 2]{Section 2} is devoted to deriving character identities and combinatorial explanations. In \hyperref[Section 3]{Section 3} we communicate spaces $V_{(i), \sqrt{N}}$ and $\Fib_{\infty}^{(\theta,\, l)}$, therefore interpreting obtained identities as character identities for $V_{(i), \sqrt{N}}$. 
\\\\
The author is grateful to
\begin{itemize}
\item Prof. Boris Feigin for the problem statement and valuable discussions,
\item Prof. Ole Warnaar for the valuable comments on the draft of this article.
\end{itemize}

\section{Fibonacci configurations}
\label{Section 1}
\subsection{Fibonacci-$1$ configurations}
\begin{figure}[h!]
\begin{center}
\includegraphics[scale = 0.8]{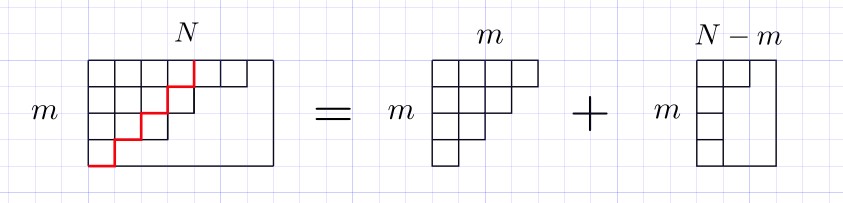}
\caption{Illustration of \eqref{qBinom1}: any Young diagram corresponding to the partition into m distinct $\leq N$ parts can be obtained by appending Young diagram contained in the $m\times (N - m)$ rectangle to the $(1, 2, \ldots, m)$ shape  from the right.}
\label{Dissection_Diagram_1}
\end{center}
\end{figure}

\begin{figure}[h!]
\begin{center}
\includegraphics[scale = 0.7]{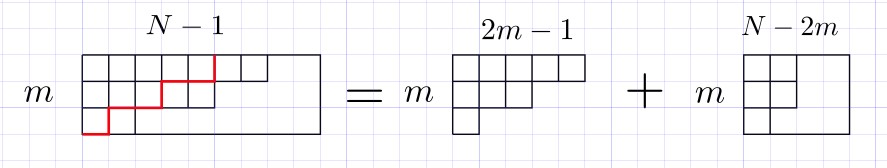}
\caption{Illustration of \eqref{qBinom2}: any Young diagram corresponding to the partition into $m$ distinct~$\leq~N~-~1$ parts with adjacent differing $\geq 2$ can be obtained by appending Young diagram contained in the \\$m\times (N - 2m)$ rectangle to the $(1, 3, \ldots, 2m - 1)$ shape  from the right.}
\label{Dissection_Diagram_2}
\end{center}
\end{figure}

\begin{definition}
Finite Fibonacci-$1$ configuration on $n$ points is a function \\$a\colon\{0, 1, \ldots, n - 1\}\longrightarrow \{0, 1\}$, such that $a_i + a_{i + 1} \leq 1$. The set of all finite Fibonacci-$1$ configurations on $n$ points is denoted by $\Fib_n^1$.  
\end{definition}
The character of the space of Fibonacci-$1$ configurations on $n$ points is given by the polynomial
\ben
\ch \Fib_n^1 := \chi_n^{1} (z, q) = \sum\limits_{a\in \Fib_n^1}\,z^{a_0 + \dots + a_{n - 1}}\,q^{0 a_0 + 1 a_1 + 2 a_2 + \dots + (n - 1) a_{n - 1}}.
\een
Fibonacci-$1$ configuration may be viewed as a set of particles placed on $n$ points in such a way that the distance between any two particles (numbers $i_1$ and $i_2$ s. t. $a_{i_1} = a_{i_2} = 1$) is strictly greater than 1, i.e.\ $|i_1 - i_2| > 1$. Then $z^{a_0 + \ldots + a_{n - 1}}$ can be interpreted as a charge of configuration and $q^{0\cdot a_0 + 1\cdot a_1 + 2\cdot a_2 + \ldots (n - 1)\cdot a_{n - 1}}$ as an energy.  
\\
There is an evident recurrent relation
\ben
\chi_{n + 1}^1 (z, q) = \chi_{n}^1(z, q) + \chi_{n - 1}^1(z, q)\cdot z\,q^n,
\een
which is based on the consideration of configurations with $a_n = 0$ and $a_n = 1$ separately.
\\
Number of Fibonacci-$1$ configurations is related to the character as
\ben
|\text{Fib}_n^1| = \chi_n^1(1, 1).
\een
Recurrent relation for this sequence is
\ben
\chi_{n + 1}^1(1, 1) = \chi_n^1(1, 1) + \chi_{n - 1}^1(1, 1),
\een
and $\chi_1(1, 1) = 2,\, \chi_2(1, 1) = 3$, i.e.\ sequence of $\chi_n^1(1, 1)$ is (up to a shift of index) a Fibonacci sequence. This fact justifies the name. 
Several first characters are
\ben
\begin{aligned}
&\chi_1^1(z, q) = 1 + z,
\\
&\chi_2^1(z, q) = 1 + z(1 + q),
\\
&\chi_3^1(z, q) = 1 + z(1 + q + q^2) + z^2\,q^2,
\\
&\chi_4^1(z, q) = 1 + z(1 + q + q^2 + q^3 + q^4) + z^2\,q^2(1 + q + q^2),
\end{aligned}
\een
We are interested in explicit formula for $\chi_{n}^{1}(z, q)$.
First of all, recall Newton's $q$-binomial theorem:
\beq
\label{qBinom1}
\prod\limits_{j = 1}^{N}\,(1 + z\,q^j) = \sum\limits_{n = 0}^{\infty}\,\sum\limits_{m = 0}^{\infty}\,p(n \,|\, m \text{ distinct parts, each } \leq N)\,z^m\,q^n = \sum\limits_{m = 0}^{N}\,z^m\,q^{\frac{m(m + 1)}{2}}\,\qbinom{N}{m}_q,
\ee
where 
\ben
\qbinom{N}{m}_q = \frac{(q)_{N}}{(q)_m\,(q)_{N - m}}, \hspace{5mm} (q)_n = \prod\limits_{i = 1}^{n} (1 - q^i) \hspace{5mm} n\in\mathbb{N}_0 \cup \{\infty\},
\een
is $q$-binomial coefficient (see \cite{AndrewsBook} for details).
This identity is explained combinatorially by Figure~\ref{Dissection_Diagram_1}.
\\
\\
Using this observation 
\beq
\label{qBinom2}
\chi_{N - 1}^{1} (zq, q)
= \sum\limits_{n = 0}^{\infty}\,\sum\limits_{m = 0}^{\infty}\,p(n \,|\, m \text{ distinct parts, each } \leq N - 1, \text{ adjacent parts differs } \geq 2)\,z^m\,q^n.
\ee
Coefficient on $z^m$ term in this sum equals
\ben
q^{1 + 3 + \ldots + (2m - 1)}\,\qbinom{N - m}{m}_q = q^{m^2}\, \qbinom{N - m}{m}_q,
\een
which is clear by ``cutting" the Young diagram as illustrated in Figure \ref{Dissection_Diagram_2}.
Thus,
\beq
\label{CharFib1}
\chi^{1}_{N - 1}\,(z\,q, q) = \sum\limits_{m = 0}^{[\frac{N}{2}]}\,z^m\,q^{m^2}\,\qbinom{N - m}{m}_q.
\ee
Formula \eqref{CharFib1} is the $q$-analog of the well-known identity
\ben
F_{n + 1} = \sum\limits_{k = 0}^{[\frac{n}{2}]} \binom{n - k}{k},
\een
\noindent
where $F_n$ is the $n$th Fibonacci number. The expression \eqref{CharFib1} also naturally arises as the generating function of lattice paths \cite{Bressoud1989} and the generating function of cylindric partitions \cite{Warnaar2023}.	

\subsection{Fibonacci-$l$ configurations}

\begin{definition}
Let $l$ be a non-negative integer. Then a finite Fibonacci-$l$ configuration on $n$ points is a function $a\colon\{0, 1, \ldots, n - 1\}\longrightarrow \{0, 1\}$, such that $a_i + a_j \leq 1$ for all $|i - j|\leq l$. The set of all Fibonacci-$l$ configurations on $n$ points is denoted by $\Fib^{l}_n$.  
\end{definition}
\noindent
The character of the space of Fibonacci-$l$ configurations on $n$ points is given by the polynomial
\ben
\ch \Fib_n^l := \chi_n^{l} (z, q) = \sum\limits_{a\in \Fib_n^l}\,z^{a_0 + \dots + a_{n - 1}}\,q^{0 a_0 + 1 a_1 + 2 a_2 + \dots + (n - 1) a_{n - 1}}.
\een
Corresponding recurrent relation is
\ben
\chi_{n + 1}^l(z, q) = \chi_{n}^l(z, q) + z\,q^n\chi_{n - l}^l(z, q) 
\een
Repeating the argument with cutting Young diagrams, we may deduce that
\ben
\chi_{n}^{l}(z, q) = P^l_{n + l}(z, q),
\een
where
\ben
P^l_{n}(z, q) = \sum\limits_{m = 0}^{\left[\frac{n}{l + 1}\right]} z^m\,q^{\frac{l + 1}{2}m(m - 1)} \qbinom{n - lm}{m}_q.
\een

\subsection{Infinite Fibonacci configurations}
\begin{definition}
Infinite Fibonacci configuration of type $(\theta,\, l)$, $\theta, l\in\mathbb{Z} $, $l \geq\theta \geq 0$, is a function $a\colon \mathbb{Z}\longrightarrow \{0, 1\}$, such that:
\begin{enumerate}
\item $a_i + a_{i + 1} + a_{i + 2} + \ldots + a_{i + l} \leq 1$,
\item $a_i = 0$\hspace{2mm} for $i \gg 1$,
\item $a_{-n} = 1$, if $-n~\equiv~\theta \mod l + 1$, $a_{-n} = 0$ otherwise for $n \gg 1$.
\end{enumerate}
The set of all infinite Fibonacci configurations of type $(\theta,\, l)$ is denoted as $\Fib^{(\theta,\,l)}_{\infty}$. The configuration which takes the value 1 only on elements of type $\theta - (l + 1)k$, $k > 0$ is called the vacuum configuration and denoted as $\Vac_{(\theta,\,l)}$. Notice that $\Fib^{(0, 0)}_{\infty}$ may be naturally identified with a basis of Fermionic Fock space (space of semi-infinite wedge forms).
\end{definition}
The character of the space of such configurations is defined as
\ben
\ch\Fib^{(\theta,\, l)}_{\infty} := \chi_{\infty}^{(\theta,\,l)}(z, q) = \sum\limits_{\substack{a\in \Fib^{(\theta,\,l)}_{\infty}\\b := a - \Vac_{(\theta,\,l)}}} \prod\limits_{i\in\mathbb{Z}} z^{b_i}\, q^{i b_i},
\een
where $a - \Vac_{(\theta,\,l)}$ means difference of the two functions on $\mathbb{Z}$. Therefore, elements of $\Fib^{(\theta,\,l)}_{\infty}$ may be viewed as configurations of particles on $\mathbb{Z}$ with the vacuum state $\Vac_{(\theta,\,l)}$. Every particle has charge 1 and energy equal to its position.
Consider limit
\ben
\lim\limits_{m \to +\infty} P^l_m(z, q) = \sum\limits_{n = 0}^{\infty} \frac{z^n q^{(l + 1)\frac{n(n - 1)}{2}}}{(q)_n} =: P^l_{\infty}(z, q).
\een
Denote
\ben
P^{(\theta,\,l)}_{\infty}(z, q) := P^l_{\infty}(zq^{\theta}, q).
\een
Note that $P^1_{\infty}(zq, q)$ and $P^1_{\infty}(zq^2, q)$ are the left sides of the famous Rogers-Ramanujan identities \cite{Rogers1894, RogRam, Schur}.  $P^{(\theta,\,l)}_{\infty}(z, q)$ is the character of subset $W^{(\theta,\,l)}_{\theta - (l + 1)}\subset\Fib^{(\theta,\,l)}_{\infty}$, consisting of the configurations stabilized on number $\theta - (l + 1)$. By the configuration stabilized on number $\theta - (l + 1)$ we mean an element of $\Fib^{(\theta,\,l)}_{\infty}$ which takes value 1 on all integer numbers of type $\theta - m(l + 1) \leq \theta - (l + 1)$, $m\in\mathbb{N}$ and value 0 on numbers  not equal to $\theta - m(l + 1)$, but $\leq \theta - (l + 1)$, $m\in\mathbb{N}$.
\\ 
The character of configurations, consisting of elements, stabilizing on $\theta - k(l + 1)$ is
\ben
\begin{aligned}
\chi_{\infty}^{(\theta,\,l)}|_{W^{(\theta,\,l)}_{\theta - k(l + 1)}}(z, q) &= z^{k - 1}\,q^{\sum\limits_{i = 1}^{k - 1} \theta - i(l + 1)} P^{(\theta,\,l)}_{\infty}(z\,q^{-(l + 1)(k - 1)}, q) 
\\
&= \sum_{m = -k + 1}^{\infty} \frac{(zq^{\theta})^m\, q^{(l + 1)\frac{m(m - 1)}{2}}}{(q)_{m + k - 1}}, 
\end{aligned}
\een
Then whole space's character may be calculated via the following limit procedure
\ben
\chi_{\infty}^{(\theta,\,l)}(z, q) = \lim\limits_{k \to +\infty} \chi_{\infty}^{(\theta,\,l)}|_{W^{(\theta,\,l)}_{\theta - k(l + 1)}} (z, q) = \sum\limits_{m\in\mathbb{Z}} \frac{(zq^{\theta})^m\, q^{(l + 1)\frac{m(m - 1)}{2}}}{(q)_{\infty}}.
\een
\section{Character identities}
\label{Section 2}
General idea of this section is to calculate the character of $\Fib_{\infty}^{(\theta,\,l)}$ by cutting $\mathbb{Z}$ into three pieces:
\ben
\mathbb{Z} = \mathbb{Z}_{< 0}\sqcup\{0, 1, \ldots, l - 1\}\sqcup\mathbb{Z}_{\geq l}.
\een  
We calculate the contribution of every part of this ``dissection" to the character and get a non-trivial character formula. Careful consideration of this formula shows its combinatorial sense~---~it is a generating function of partitions based on Durfee rectangle with the ratio of the sides $1:(l + 1)$. These identities were firstly studied by G. Andrews in \cite{Andrews_Part1}. The suggested procedure may be viewed as an alternative proof of these identities. Moreover, derived identities generalize identities of G. Andrews, as soon as they admit constant shift in base and height of Durfee rectangles.
\\\\
Firstly, the simplest cases $l = 0$ and $l = 1$ are considered. After that general case $l = k\in\mathbb{Z_+}$ is examined.
\subsection{$l = 0$ case}
In this case the distance condition is trivial, therefore, the character of the whole space $\Fib_{\infty}^{(0, 0)}$ is product of the left and right part characters. 
The character of the right part (corresponding to the particles on $\mathbb{Z}_+$) is
\ben
\label{CharFerm}
\lim\limits_{n\to\infty} \prod\limits_{j = 0}^{n} (1 + zq^j) = \lim\limits_{n\to\infty} \sum\limits_{m = 0}^{n + 1} z^m q^{\frac{m(m - 1)}{2}}\qbinom{n + 1}{m}_q = \sum\limits_{m = 0}^{\infty} \frac{z^m q^{\frac{m(m - 1)}{2}}}{(q)_m}.
\een
The character of the left part (corresponding to the particles on $\mathbb{Z}\setminus\mathbb{Z}_+$) is found by the same way:
\ben
\lim\limits_{n\to\infty} \prod\limits_{j = 1}^{n} (1 + z^{-1}q^j) = \lim\limits_{n\to\infty} \sum\limits_{k = 0}^{n} z^{-l} q^{\frac{k(k + 1)}{2}}\qbinom{n}{k}_q = \sum\limits_{l = 0}^{\infty} \frac{z^{-k} q^{\frac{k(k + 1)}{2}}}{(q)_k}.
\een
Finally, we get a nontrivial identity (equivalent to the well-known Jacobi's triple product identity \cite{Jacobi_Fundamenta})
\ben
\chi^{(0, 0)}_{\infty} = \sum\limits_{n \in \mathbb{Z}} \frac{z^n\,q^{\frac{n(n - 1)}{2}}}{(q)_{\infty}} = \sum\limits_{k = 0}^{\infty} \frac{z^{-k} q^{\frac{k(k + 1)}{2}}}{(q)_k} \sum\limits_{m = 0}^{\infty} \frac{z^m q^{\frac{m(m - 1)}{2}}}{(q)_m}, 
\een
which could be rewritten by taking the coefficient of $z^s$ on both sides:
\beq
\label{DurfeeRectangleIdentity0}
z^s: \hspace{5mm}\frac{1}{(q)_{\infty}} = \sum\limits_{\substack{m, k \geq 0 \\ m - k = s}} \frac{q^{mk}}{(q)_m (q)_k},
\ee
where $s$ is a fixed number.

\begin{figure}[h!]
\begin{center}
\includegraphics[scale = 1]{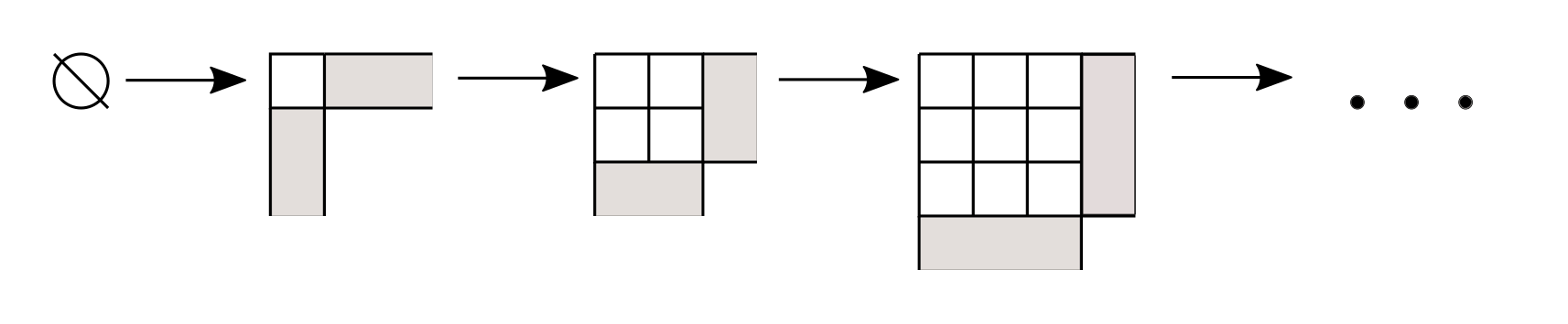}
\caption{Identity \eqref{DurfeeRectangleIdentity0} at $s = 0$. \\ Durfee rectangles $k\times k$}
\label{DurfeeSquare}
\end{center}
\end{figure}

\begin{figure}[h!]
\begin{center}
\includegraphics[scale = 1]{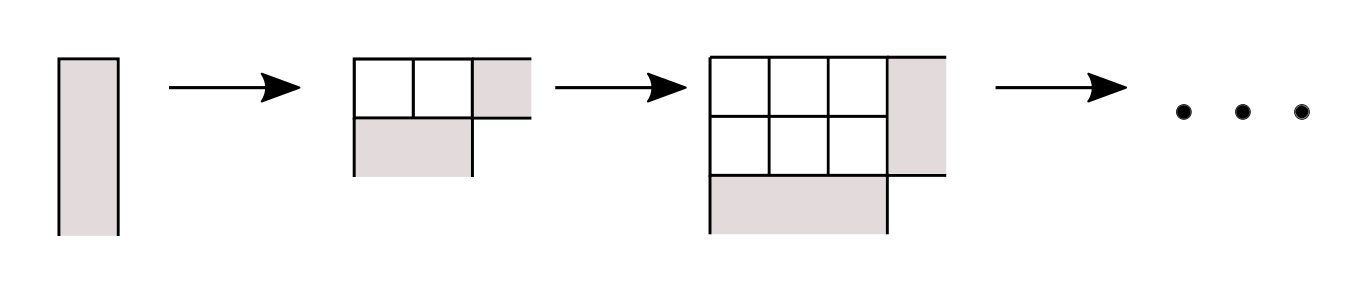}
\caption{Identity \eqref{DurfeeRectangleIdentity0} at $s = 1$. \\ Durfee rectangles $k\times (k + 1)$}
\label{DurfeeRectangle01}
\end{center}
\end{figure}

\begin{figure}[h!]
\begin{center}
\includegraphics[scale = 1]{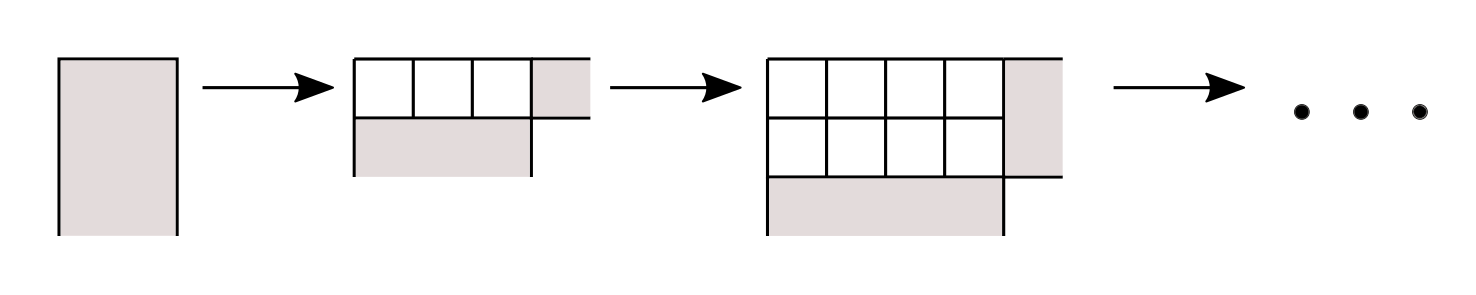}
\caption{Identity \eqref{DurfeeRectangleIdentity0} at $s = 2$. \\ Durfee rectangles $k\times (k + 2)$}
\label{DurfeeRectangle02}
\end{center}
\end{figure}
For any $s \in \mathbb{N}_0$, for any partition there is a maximal $n\in\mathbb{N}_0$, such that this partition contains rectangle $n\times(n + s)$. Denote the set of all such partitions $P_{1, s, n}$. Then the set of all partitions is 
\ben
P = \bigsqcup\limits_{n \in \mathbb{N}_0} P_{1, s, n}
\een
and there is a natural linear order on set of $P_{1, s, n}$ inherited from the linear order on $\mathbb{N}_0$. The generating function of $P_{1, s, n}$ is $\frac{q^{n(n + s)}}{(q)_{n}(q)_{n + s}}$, which gives identity \eqref{DurfeeRectangleIdentity0}. It is illustrated in particular cases in Figures \ref{DurfeeSquare}, \ref{DurfeeRectangle01} and \ref{DurfeeRectangle02}.       

\subsection{$l = 1$ case}
\noindent
In this case the distance between two particles is strictly greater than 1. Without loss of generality, consider case $(1, 1)$.  
All configurations in $\Fib_{\infty}^{(1, 1)}$ can be divided in two sets:
\ben
\Fib\limits_{\infty}^{(1, 1)} = A^{(1, 1)}_{\varnothing} \sqcup A^{(1, 1)}_{0},
\een
where $A^{(1, 1)}_{\varnothing} = \{a\in\Fib_{\infty}^{(1, 1)} |\, a_0 = 0\}$ and $A^{(1, 1)}_{0} = \{a\in\Fib_{\infty}^{(1, 1)} |\, a_0 = 1\}$.
\\
The character of $\Fib_{\infty}^{(1, 1)}$ is
\ben
\ch \Fib\limits_{\infty}^{(1, 1)} = \ch A^{(1, 1)}_{\varnothing} + \ch A^{(1, 1)}_{0}.
\een
Every configuration in $A^{(1, 1)}_{\varnothing}$ and $A^{(1, 1)}_{0}$ may be ``cut" into left and right part:
\beq
\label{LeftRightParts}
\begin{aligned}
&L^{(1, 1)}_{\varnothing} = \{j\colon \mathbb{Z}_{< 0} \rightarrow \{0, 1\}|\hspace{2mm} j = a|_{\mathbb{Z}_{< 0}},\hspace{2mm} a \in A^{(1, 1)}_{\varnothing}\},
\\
&R^{(1, 1)}_{\varnothing} = \{r\colon \mathbb{Z}_{> 0} \rightarrow \{0, 1\}|\hspace{2mm} r = a|_{\mathbb{Z}_{> 0}},\hspace{2mm} a \in A^{(1, 1)}_{\varnothing}\},
\\
&L^{(1, 1)}_{0} = \{j\colon \mathbb{Z}_{< 0} \rightarrow \{0, 1\}|\hspace{2mm} j = a|_{\mathbb{Z}_{< 0}},\hspace{2mm} a \in A^{(1, 1)}_{0}\},
\\
&R^{(1, 1)}_{0} = \{r\colon \mathbb{Z}_{> 0} \rightarrow \{0, 1\}|\hspace{2mm} r = a|_{\mathbb{Z}_{> 0}},\hspace{2mm} a \in A^{(1, 1)}_{0}\}.
\end{aligned}
\ee
Corresponding characters of these parts are
\beq
\label{LR_Char}
\begin{aligned}
&\ch L^{(1, 1)}_{\varnothing} = \sum\limits_{\substack{j\in L^{(1, 1)}_{\varnothing}\\b := j - \Vac_{(1, 1)}|_{\mathbb{Z}_{< 0}}}} \prod\limits_{i < 0}\, z^{b_i}\, q^{i b_i},
\\
&\ch R^{(1, 1)}_{\varnothing} = \sum\limits_{r\in R^{(1, 1)}_{\varnothing}} \prod\limits_{i > 0}\, z^{r_i}\, q^{i r_i},
\\
&\ch L^{(1, 1)}_{0} = \sum\limits_{\substack{j\in L^{(1, 1)}_{0}\\b := j - \Vac_{(1, 1)}|_{\mathbb{Z}_{< 0}}}} \prod\limits_{i < 0}\, z^{b_i}\, q^{i b_i},
\\
&\ch R^{(1, 1)}_{0} = \sum\limits_{r\in R^{(1, 1)}_{0}} \prod\limits_{i > 0}\, z^{r_i}\, q^{i r_i}.
\end{aligned}
\ee
Then the character of $\Fib_{\infty}^{(1, 1)}$ may be rewritten as
\beq
\label{01Identity}
\ch \Fib_{\infty}^{(1, 1)}  = \sum\limits_{n\in\mathbb{Z}} \frac{z^n\,q^{n^2}}{(q)_{\infty}} = \ch L^{(1, 1)}_{\varnothing}\ch R^{(1, 1)}_{\varnothing} + z \ch L^{(1, 1)}_{0}\ch R^{(1, 1)}_{0}
\ee
To get the character identity we need to calculate \eqref{LR_Char} explicitly. Characters of the right parts are evident:
\beq
\begin{aligned}
&\ch R^{(1, 1)}_{\varnothing} = P^{(1, 1)}_{\infty} (zq, q) = \sum\limits_{n = 0}^{\infty} \frac{z^n\,q^{n^2}}{(q)_n},  
\\
&\ch R^{(1, 1)}_{0} = P^{(1, 1)}_{\infty} (zq^2, q) = \sum\limits_{n = 0}^{\infty} \frac{z^n\,q^{n(n + 1)}}{(q)_n}.
\end{aligned}
\ee
Characters of the left parts may be calculated as
\beq
\label{LeftPart}
\begin{aligned}
&\ch L^{(1, 1)}_{\varnothing} = \lim\limits_{n\to +\infty} \frac{q^{n^2}}{z^n} \chi_{2n - 1}^{1}(zq^{-1}, q^{-1}),
\\	
&\ch L^{(1, 1)}_{0} = \lim\limits_{n\to +\infty} \frac{q^{(n + 1)^2}}{z^{n + 1}} \chi_{2n}^{1}(zq^{-2}, q^{-1})
\end{aligned}
\ee
Formula \eqref{LeftPart} deserves particular explanation. Let's consider the case of $\ch L^{(1, 1)}_{\varnothing}$. Firstly,  we  compute the character of $L^{(1, 1)}_{\varnothing}$ subspace, consisting of the configurations, stabilized on number $-(2n + 1)$ (i.e.\ $a_i = 1$ for all odd $i\leq -(2n + 1)$ and $a_i = 0$ for all even $i\leq -(2n + 1)$ with mentioned combinatorial property of the distance between 1's being $> 1$) and tend $n$ to $+\infty$. This subspace is of the configurations on $2n - 1$ points, as soon as positions are negative we put $q \rightarrow q^{-1}$, as soon as they start with $-1$ we put $z \rightarrow z\,q^{-1}$ in $\chi^{1}_{2n - 1}(z, q)$. Due to the vacuum, empty configuration has weight
\beq
z^{-n}\,q^{\sum\limits_{i = 1}^{n} (2i - 1)} = z^{-n}\,q^{n^2}, 
\ee
which explains the factor in formula \eqref{LeftPart}. All characters of the left parts in this chapter are calculated via the same procedure.
\noindent
Calculation of the $\ch L^{(1, 1)}_{\varnothing}$ gives:
\beq
\begin{aligned}
    &\frac{q^{n^2}}{z^n} \chi_{2n - 1}^{1}(zq^{-1}, q^{-1}) = \frac{q^{n^2}}{z^n} P_{2n}^1(zq^{-1}, q^{-1}) = \frac{q^{n^2}}{z^n} \sum\limits_{i = 0}^n (z\,q^{-1})^i q^{-i(i - 1)} \qbinom{2n - i}{i}_{q^{-1}} 
\\
    & = \sum\limits_{k = 0}^{n} z^{-k}\,q^{k^2}\qbinom{n + k}{2k}_q \overset{n \to +\infty}{\longrightarrow} \sum\limits_{k = 0}^{\infty} \frac{z^{-k}\,q^{k^2}}{(q)_{2k}},
\end{aligned}
\ee
where we use
\beq
\qbinom{a}{b}_{q^{-1}} = q^{-(a - b) b} \qbinom{a}{b}_q.
\ee
Similar calculation of the $\ch L^{(1, 1)}_{0}$ gives:
\beq
\begin{aligned}
    &\frac{q^{(n + 1)^2}}{z^{n + 1}} \chi_{2n}^{1}(zq^{-1}, q^{-1}) = \frac{q^{n^2}}{z^n} P_{2n + 1}^1(zq^{-2}, q^{-1}) = \frac{q^{n^2}}{z^n} \sum\limits_{i = 0}^n (z\,q^{-2})^i q^{-i(i - 1)} \qbinom{2n + 1 - i}{i}_{q^{-1}}
    \\
    & = \frac{1}{z}\,\sum\limits_{k = 0}^{n} z^{-k}\,q^{(k + 1)^2}\qbinom{n + k + 1}{n - k}_q \overset{n \to \infty}{\longrightarrow} \frac{1}{z}\,\sum\limits_{k = 0}^{\infty} \frac{z^{-k}\,q^{(k + 1)^2}}{(q)_{2k + 1}}.
\end{aligned}
\ee
Substituting these results in \eqref{01Identity}, we get an identity in the formal power series ring $\mathbb{Z}[z^{\pm 1}][[q]]$:
\beq
\label{01IdentityExplicit}
\sum\limits_{n\in\mathbb{Z}} \frac{z^n\,q^{n^2}}{(q)_{\infty}} = \sum\limits_{n = 0}^{\infty} \frac{z^n\,q^{n^2}}{(q)_n} \sum\limits_{k = 0}^{\infty} \frac{z^{-k}\,q^{k^2}}{(q)_{2k}} + \sum\limits_{n = 0}^{\infty} \frac{z^n\,q^{n(n + 1)}}{(q)_n} \sum\limits_{k = 0}^{\infty} \frac{z^{-k}\,q^{(k + 1)^2}}{(q)_{2k + 1}}.
\ee
Similar calculations for $(0, 1)$ case give:
\beq
\label{11IdentityExplicit}
\sum\limits_{n\in\mathbb{Z}} \frac{z^n\,q^{n(n - 1)}}{(q)_{\infty}} = \sum\limits_{n = 0}^{\infty} \frac{z^{n + 1}\,q^{n(n + 1)}}{(q)_n} \sum\limits_{k = 0}^{\infty} \frac{z^{-k}\,q^{k(k + 1)}}{(q)_{2k}} +
\sum\limits_{n = 0}^{\infty} \frac{z^n\,q^{n^2}}{(q)_n} \sum\limits_{k = 0}^{\infty} \frac{z^{-k}\,q^{k(k + 1)}}{(q)_{2k + 1}}.
\ee
Equating coefficients on fixed degrees of $z$ variable, identity \eqref{01IdentityExplicit} may be rewritten as series of identities:

\beq
\label{11pos}
z^s: \hspace{5mm}\frac{1}{(q)_{\infty}} = \sum\limits_{k = 0}^{\infty} \frac{q^{2k(k + s)}}{(q)_{k + s}\,(q)_{2k}} + \sum\limits_{k = 0}^{\infty} \frac{q^{(2k + 1)(k + s + 1)}}{(q)_{k + s}\,(q)_{2k + 1}}  = \sum\limits_{k = 0}^{\infty} \frac{q^{2k(k + s)}(1 + q^{3k + s + 1} - q^{2k + 1})}{(q)_{k + s}(q)_{2k + 1}} \hspace{5mm} s \in \mathbb{N}_0,
\ee
\beq
\label{11neg}
\begin{aligned}
z^{-s}: \hspace{5mm}\frac{1}{(q)_{\infty}} &= \sum\limits_{n = 0}^{\infty} \frac{q^{2n(n + s)}}{(q)_{n}\,(q)_{2(n + s)}} + \sum\limits_{n = 0}^{\infty} \frac{q^{(2(n + s) + 1)(n + 1)}}{(q)_{n}\,(q)_{2(n + s) + 1}} 
\\ 
&= \sum\limits_{n = 0}^{\infty} \frac{q^{2n(n + s)}(1 - q^{2(n + s) + 1} + q^{2(n + s) + n + 1})}{(q)_{n}(q)_{2(n + s) + 1}} \hspace{5mm}  s \in \mathbb{N}_0.
\end{aligned}
\ee
Identity \eqref{11IdentityExplicit} may be analogously rewritten as series of identities:
\beq
\label{01pos}
\begin{aligned}
z^s: \hspace{5mm}\frac{1}{(q)_{\infty}} &= \frac{1}{(q)_{s - 1}} + \sum\limits_{k = 0}^{\infty} \frac{q^{(2k + 1)(k + s)}}{(q)_{2k + 1}(q)_{k + s}} + \sum\limits_{k = 0}^{\infty}\frac{q^{(2k + 2)(k + s + 1)}}{(q)_{2k + 2}(q)_{k + s}}
\\
&= \sum\limits_{k = 0}^{\infty} \frac{q^{(2k + 1)(k + s)}(1 - q^{2k + 2} + q^{3k + s + 2})}{(q)_{k + s} (q)_{2k + 2}}  \hspace{5mm} s \in \mathbb{N},
\end{aligned}
\ee
\beq
\label{01neg}
\begin{aligned}
z^{-s}: \hspace{5mm}\frac{1}{(q)_{\infty}} &= \sum\limits_{n = 0}^{\infty} \frac{q^{2(n + 1)(n + s + 1)}}{(q)_{n}\,(q)_{2(n + s + 1)}} + \sum\limits_{n = 0}^{\infty} \frac{q^{(2(n + s) + 1)n}}{(q)_{n}\,(q)_{2(n + s) + 1}}
\\ 
&=\sum\limits_{n = 0}^{\infty} \frac{q^{2n(n + s) + n}(1 + q^{2(n + s + 1) + n} - q^{2(n + s + 1)})}{(q)_{n}(q)_{2(n + s + 1)}} \hspace{5mm} &s \in \mathbb{N}_0.
\end{aligned}
\ee
\\
\begin{remark}
Using convention $\frac{1}{(q)_n} = 0$ if $n < 0$ \eqref{11pos} and \eqref{11neg} may be rewritten as \eqref{11pos} for all $s\in\mathbb{Z}$. The same argument applies to \eqref{01pos} and \eqref{01neg}. The reason of writing these identities separately is more evident structure of summands and simpler classification, which will be explained later.      
\end{remark}
In the case $s = 0$ identities \eqref{11pos} and \eqref{11neg} turn into:
\beq
\label{AndrewsIdentity1}
\frac{1}{(q)_{\infty}} = \sum\limits_{k = 0}^{\infty} \frac{q^{2k^2}(1 + q^{3k + 1} - q^{2k + 1})}{(q)_k(q)_{2k + 1}} = \sum\limits_{k = 0}^{\infty} \frac{q^{2k^2}}{(q)_k (q)_{2k}} + \sum\limits_{k = 0}^{\infty} \frac{q^{(k + 1)(2k + 1)}}{(q)_k (q)_{2k + 1}},
\ee
This identity was firstly obtained by Andrews in \cite{Andrews_Part1} by considering Durfee rectangle.

\begin{figure}[!htb]
\minipage{0.48\textwidth}
  \includegraphics[width=\linewidth]{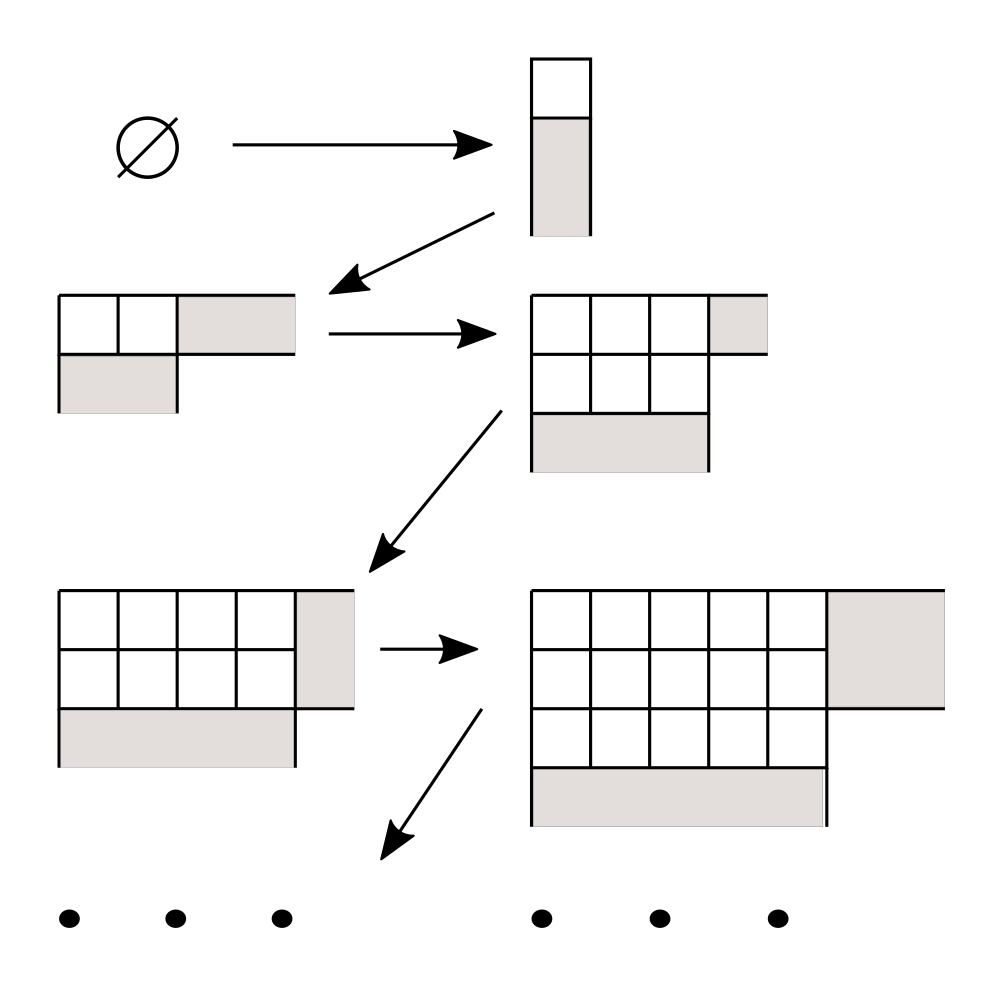}
  
  \caption{Identity \eqref{AndrewsIdentity1}\\ Durfee rectangles: $k\times 2k$ \\
  Enveloping rectangles: $(k + 1)\times(2k + 1)$}
  \label{Picture 1}
\endminipage\hfill
\minipage{0.48\textwidth}
  \includegraphics[width=\linewidth]{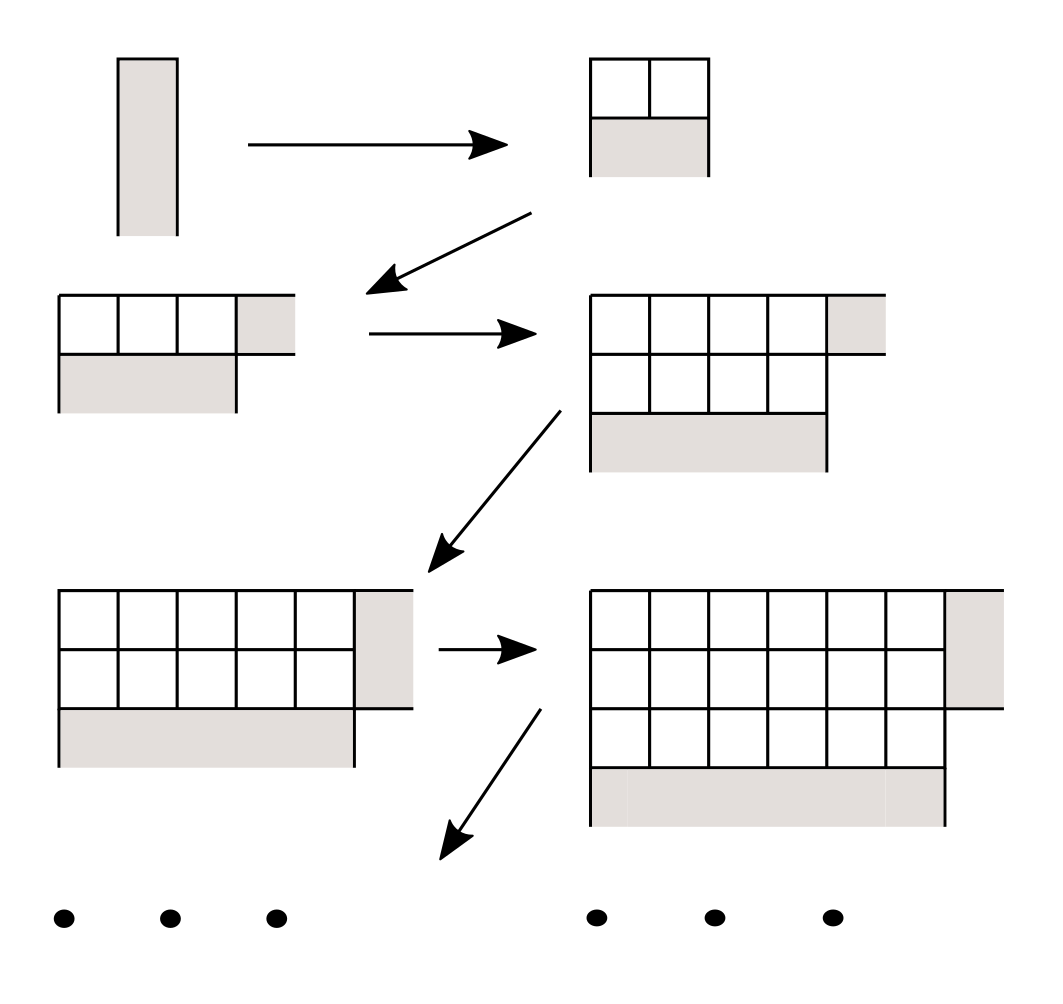}
  \caption{Identity \eqref{01neg} at s = 0\\ Durfee rectangles $k\times(2k + 1)$ \\ Enveloping rectangles: $(k + 1)\times(2k + 2)$}\label{Picture 2}
\endminipage\hfill

\end{figure}

\begin{figure}[!htb]
\minipage{0.48\textwidth}
  \includegraphics[width=\linewidth]{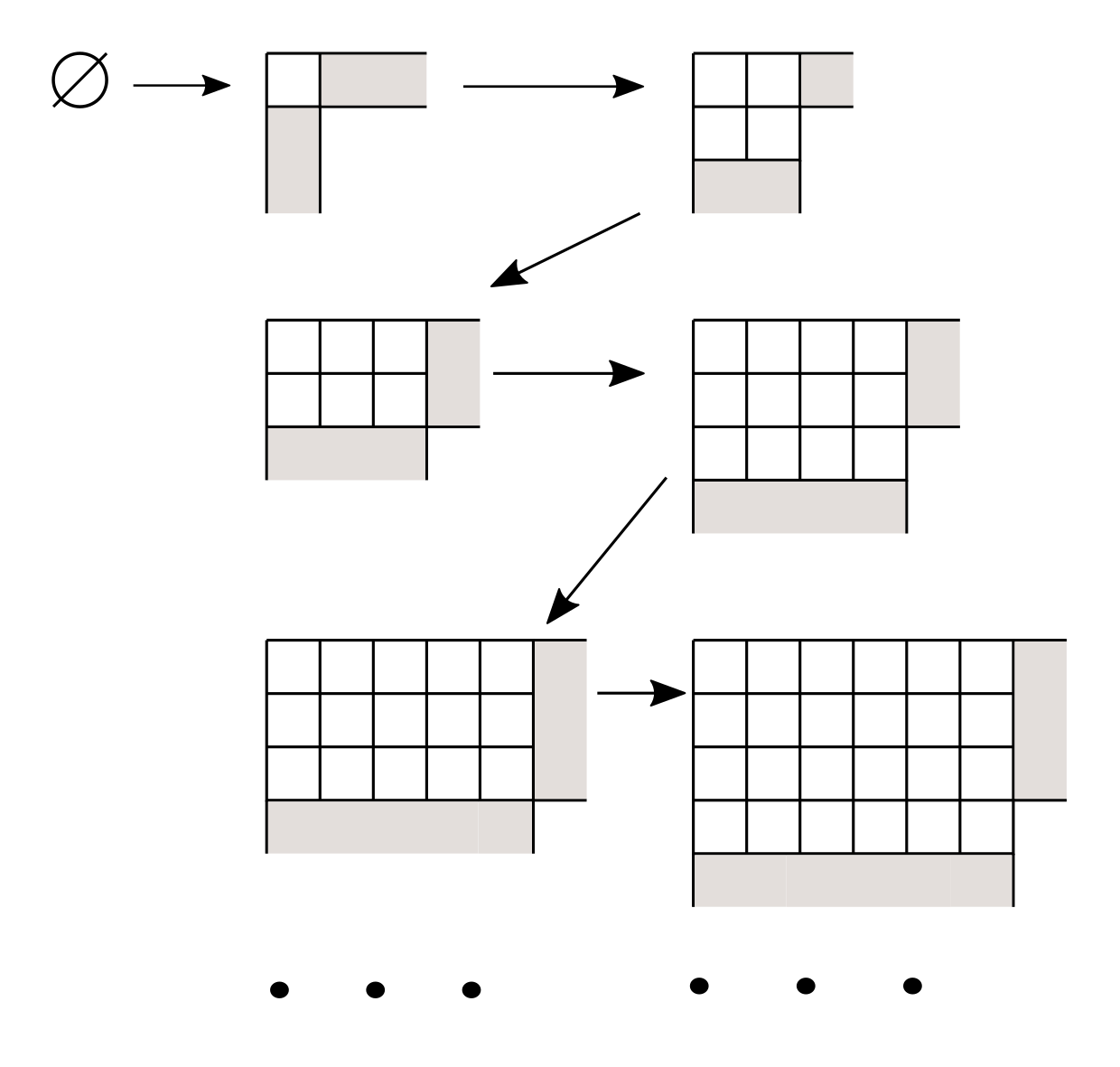}
  \caption{Identity \eqref{01pos} at s = 1 \\ Durfee rectangles $(k + 1)\times(2k 
 + 1)$ \\ Enveloping rectangles: $(k + 2)\times(2k + 2)$}\label{Picture 3}
\endminipage\hfill
\minipage{0.48\textwidth}
  \includegraphics[width=\linewidth]{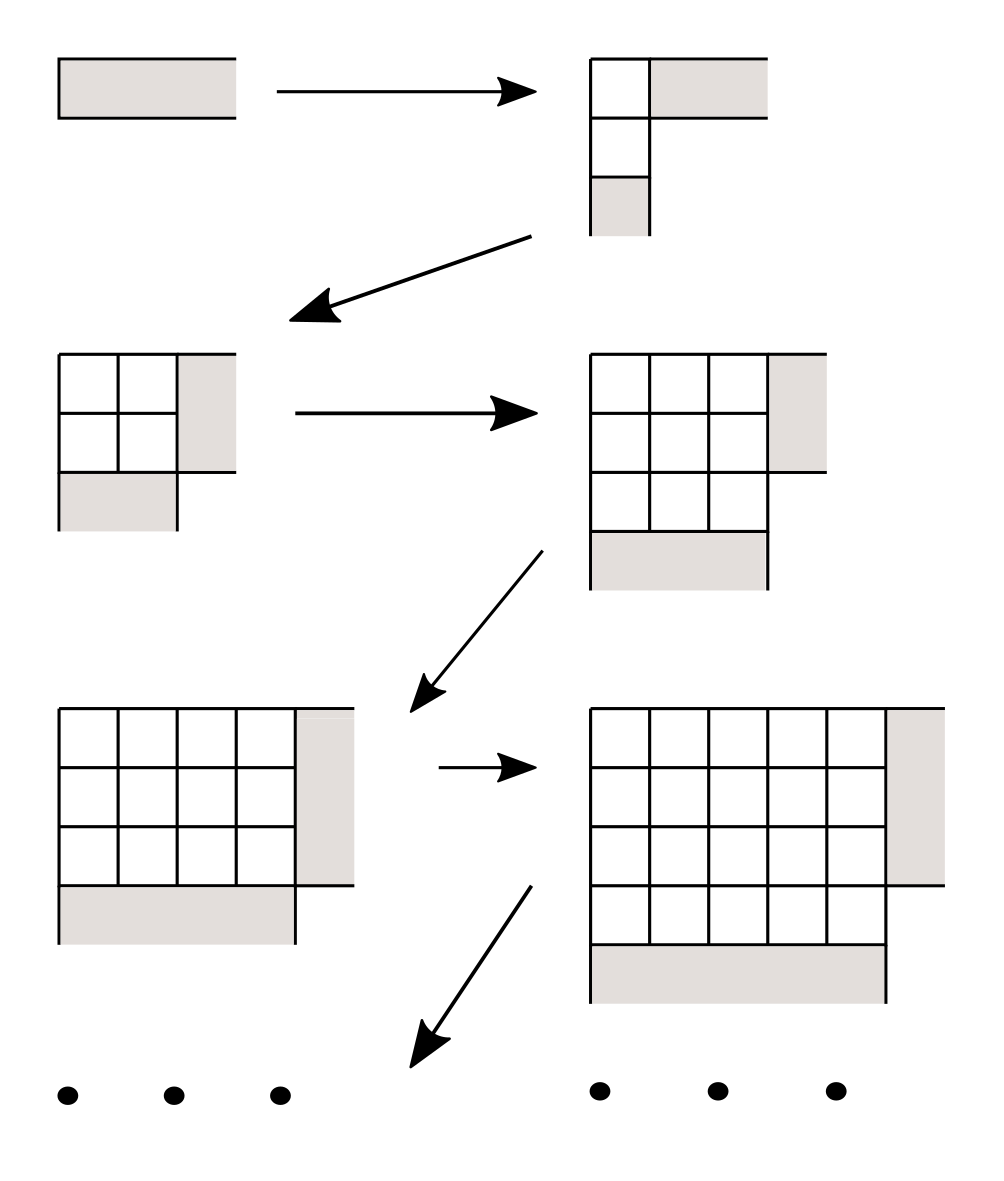}
  \caption{Identity \eqref{11pos} at s = 1 \\ Durfee rectangles: $(k + 1)\times 2k$ \\ Enveloping rectangles: $(k + 2)\times(2k + 1)$}\label{Picture 4}
\endminipage\hfill

\end{figure}

For any $n, m \in \mathbb{N}_0$, for any partition it either does not contain $n\times m$ rectangle or there is maximal $k\in\mathbb{N}_0$, such that this partition contains rectangle $(k + n)\times(2k + m)$ (such rectangle we call Durfee rectangle of this partition). Set of all partitions from the first case we denote $\tilde{P}_{2, n, m}$. The generating function of $\tilde{P}_{2, n, m}$ is
\ben
\chi_{\substack{n \geq 1 \\ m\geq 1}} \sum\limits_{j = 0}^{m - 1} \frac{q^{nj}}{(q)_{n - 1}(q)_j}.
\een
Denote the set of all partitions with Durfee rectangle $(k + n)\times(2k + m)$ as $P_{2, k, n, m}$. Then the set of all partitions is 
\ben
P = \tilde{P}_{2, n, m}\bigsqcup\left(\bigsqcup\limits_{k \in \mathbb{N}_0} P_{2, k, n, m}\right).
\een
Partitions from $P_{2, k, n, m}$ contain $(k + n)\times(2k + m)$ rectangle but do not contain $(k + 1 + n)\times(2(k + 1) + m)$ rectangle (which we call enveloping, G. Andrews in \cite{Andrews_Part1} calls them ambient). $P_{2, k, n, m}$ can be divided in two sets: $P_{2, k, n, m}^{0}$~---~partitions, which do not contain $(k + n + 1)\times(2k + 1 + m)$ rectangle, and $P_{2, k, n, m}^{1}$ containing $(k + n + 1)\times(2k + 1 + m)$ rectangle with the generating functions
\ben
\frac{q^{(k + n)(2k + m)}}{(q)_{k + n}(q)_{2k + m}} \hspace{5mm}\text{and} \hspace{5mm} \frac{q^{(k + n + 1)(2k + m + 1)}}{(q)_{k + n}(q)_{2k + m + 1}}
\een
respectively. Thus, we get general family of such type identities: 
\ben
\frac{1}{(q)_{\infty}} = \sum\limits_{k = 0}^{\infty} \frac{q^{(k + n)(2k + m)}}{(q)_{k + n} (q)_{2k + m}} + \sum\limits_{k = 0}^{\infty} \frac{q^{(k + n + 1)(2k + m + 1)}}{(q)_{k + n}(q)_{2k + m + 1}} + \chi_{\substack{n \geq 1 \\ m\geq 1}} \sum\limits_{j = 0}^{m - 1} \frac{q^{nj}}{(q)_{n - 1}(q)_j}
\een
for $n, m\in \mathbb{N}_0$. Several cases with the natural ordering on $P$'s are illustrated in Figures \ref{Picture 1}, \ref{Picture 2}, \ref{Picture 3} and \ref{Picture 4}. 

\subsection{General case}
\begin{figure}[h!]
\begin{center}
\includegraphics[scale = 0.5]{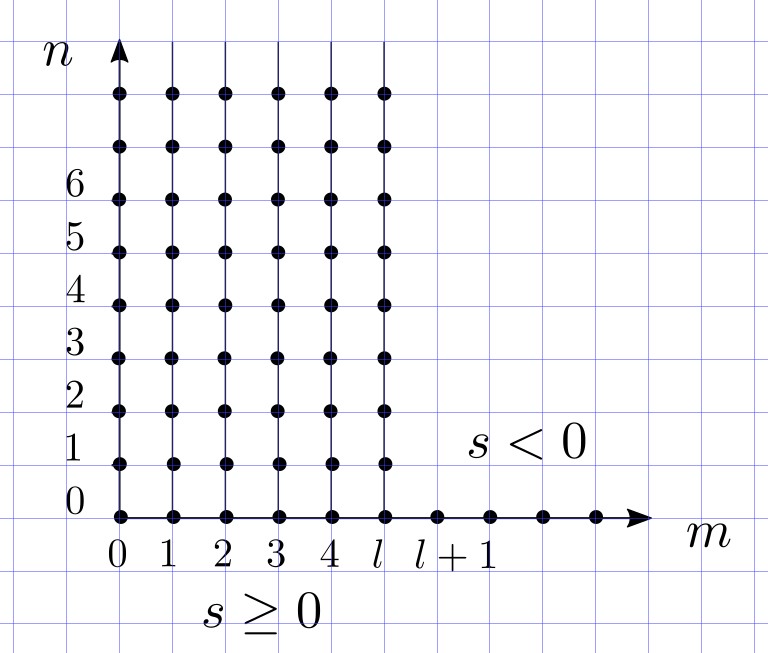}

\caption{Identities of type \eqref{Final identity l + 1} obtained as character identities for $\Fib_{\infty}^{(\theta,\, l)}$ are marked by black dots. $l$ is fixed, $0 \leq \theta \leq l$.}

\label{Diagram}
\end{center}
\end{figure}

Deriving character identities in general case is similar to the $l = 1$ case. Consider splitting $\mathbb{Z}$ into three parts 
\beq
\label{splittingZ}
\mathbb{Z} = \mathbb{Z}_{< 0}\sqcup\{0, 1, \ldots, l - 1\}\sqcup\mathbb{Z}_{\geq l}. 
\ee
Denote $A^{(\theta,\,l)}_{\varnothing} = \{a\in\Fib_{\infty}^{(\theta,\,l)} |\, a_i = 0 \, \text{ for } 0 \leq i \leq l - 1\}$, $A^{(\theta,\,l)}_{i} = \{a\in\Fib_{\infty}^{(\theta,\,l)} |\, a_i = 1\}$. Then the character of $\Fib^{(\theta,\,l)}_{\infty}$ splits into $l + 1$ summands: 

\beq
\label{CharIdThetaL}
\begin{aligned}
&\ch\Fib^{(\theta,\,l)}_{\infty} = \ch A^{(\theta,\, l)}_{\varnothing} + \ch A^{(\theta,\, l)}_{0} + \ch A^{(\theta,\, l)}_{1} + \dots + \ch A^{(\theta,\, l)}_{l - 1}
\\
&= \ch L^{(\theta,\, l)}_{\varnothing}\ch R^{(\theta,\, l)}_{\varnothing} + z\ch L^{(\theta,\, l)}_{0}\ch R^{(\theta,\, l)}_{0} + zq\ch L^{(\theta,\, l)}_{1}\ch R^{(\theta,\, l)}_{1} + \dots + zq^{l - 1}\ch L^{(\theta,\, l)}_{l - 1}\ch R^{(\theta,\, l)}_{l - 1},
\end{aligned}
\ee
where the left and right parts and their characters are defined by analogy with \eqref{LeftRightParts} and \eqref{LR_Char} with respect to splitting \eqref{splittingZ}.
\\ 
The characters of the right parts are
\ben
\ch R^{(\theta,\, l)}_{\varnothing} = P^{l}_{\infty}(z\,q^l, q), \hspace{5mm} \ch R^{(\theta,\, l)}_{i} = P^{l}_{\infty}(z\,q^{i + l + 1}, q),
\een
where
\ben
P^{l}_{\infty} (z, q) = \sum\limits_{n = 0}^{\infty} \frac{z^n\,q^{\frac{l + 1}{2}n(n - 1)}}{(q)_n}. 
\een
Denote 
\ben
\ch L^{(\theta,\, l)}_{\varnothing} = \Tilde{P}^{(\theta,\, l),\, -1}_{\infty}, \hspace{5mm} \ch L^{(\theta,\, l)}_{i} = \Tilde{P}^{(\theta,\, l), i - l - 1}_{\infty}.  
\een
The strategy is to calculate $\Tilde{P}^{(0, l), - k}_{\infty}$, where $-k$ ranges from $-1$ to $- l - 1$ and shows the maximal number where 1 may appear in the corresponding configurations on $\mathbb{Z}_{< 0}$. For simplicity, the case $\theta = 0$ is considered first. Denote $W_{k, b}^{(0, l)}\subset L_{l + 1 - k}^{(0, l)}$ ($L_{l}^{(0, l)}\equiv L_{\varnothing}^{(0, l)}$) the set of configurations stabilizing on number $-(l + 1)(b + 1)$,\, $b\in\mathbb{N}_0$. Following limit procedure gives an explicit expression for $\Tilde{P}^{(0, l), - k}_{\infty}$:
\ben
\Tilde{P}^{(0, l), - k}_{\infty} = \lim\limits_{b\to +\infty} \ch W_{k, b}^{(0, l)} = \lim\limits_{b\to +\infty} \frac{q^{\frac{(l + 1)}{2}b(b + 1)}}{z^b} \chi^{l}_{(l + 1)b - k + 1} (z\,q^{-k}, q^{-1}). 
\een
\ben
\begin{aligned}
&\frac{q^{\frac{(l + 1)}{2}b(b + 1)}}{z^b} \chi^{l}_{(l + 1)b - k + 1} (z\,q^{-k}, q^{-1}) = \frac{q^{\frac{(l + 1)}{2}b(b + 1)}}{z^b} P^{l}_{(l + 1)b - k + 1 + l} (z\,q^{-k}, q^{-1})
\\
&= \frac{q^{\frac{(l + 1)}{2}b(b + 1)}}{z^b} \sum\limits_{n = 0}^{b} z^n\,q^{-kn}\,q^{-\frac{l + 1}{2}n(n - 1)} \qbinom{(l + 1)(b + 1) - k - ln}{n}_{q^{-1}}
\\
&=\sum\limits_{m = 0}^{b} z^m\,q^{\frac{l + 1}{2}m(m + 1)}\qbinom{(l + 1)(b + 1) - k - l(b - m)}{b - m}_q \overset{b\to +\infty}{\longrightarrow} \sum\limits_{m = 0}^{\infty} \frac{z^{-m}\, q^{\frac{l + 1}{2}m(m + 1)}}{(q)_{(l + 1)(m + 1) - k}}. 
\end{aligned}
\een
In case of general $\theta$
\ben
\Tilde{P}^{(\theta,\, l), - k}_{\infty} = \lim\limits_{b\to +\infty}\, \ch W^{(\theta,\, l)}_{k, b} = \lim\limits_{b\to\infty}\, \frac{q^{\sum\limits_{k = 1}^{b} \left((l + 1)k - \theta\right)}}{z^b}\,\chi_{b(l + 1) - \theta - k + 1}\,(z\,q^{-k}, q^{-1}). 
\een
\ben
\begin{aligned}
&\chi_{b(l + 1) - \theta - k + 1}\,(z\,q^{-k}, q^{-1}) = \frac{q^{\frac{l + 1}{2}\,b(b + 1) - \theta\,b}}{z^b}\,P^l_{(b + 1)(l + 1) - \theta - k} (z\,q^{-k}, q^{-1})  
\\
& = \frac{q^{\frac{l + 1}{2}\,b(b + 1) - \theta\,b}}{z^b}\sum\limits_{n = 0}^{[b + 1 - \frac{\theta + k}{l + 1}]}\, (z\,q^{-k})^n\,q^{-\frac{l + 1}{2}\,n(n - 1)}\,\qbinom{(l + 1)(b + 1) - \theta - k - ln}{n}_{q^{-1}} 
\\
& = \sum\limits_{n = 0}^{[b + 1 - \frac{\theta + 1}{l + 1}]}\,z^{-(b - n)}\,q^{\frac{l + 1}{2}\,(b - n)(b - n + 1) - \theta(b - n)}\,\qbinom{(l + 1)(b + 1) - \theta - k - ln}{n}_{q} 
\\
& = \sum\limits_{m = \alpha(\theta, k)}^{b}\, z^{-m}\,q^{\frac{l + 1}{2}m(m + 1) - \theta\,m}\,\qbinom{(l + 1)(b + 1) - \theta - k - l(b - m)}{b - m}_{q}
\\
& = \sum\limits_{m = \alpha(\theta, k)}^{b}\, z^{-m}\,q^{\frac{l + 1}{2}m(m + 1) - \theta\,m}\,\qbinom{b + l(m + 1) + 1 - \theta - k}{(l + 1)(m + 1) - \theta - k}_{q}.
\end{aligned}
\een
Thus,
\beq
\label{PInfty}
\Tilde{P}^{(\theta,\, l), - k}_{\infty} = \sum\limits_{m = \alpha(\theta, k)}^{\infty} \frac{(z\,q^{\theta})^{-m} q^{\frac{(l + 1)}{2} m(m + 1)}}{(q)_{(m + 1)(l + 1) - \theta - k}},
\ee
where 
\ben
\begin{aligned}
&\alpha(\theta, k) = 0 \hspace{5mm} \text{if}\hspace{5mm} 1 \leq \theta + k \leq l + 1,
\\
&\alpha(\theta, k) = 1 \hspace{5mm} \text{if}\hspace{5mm}  l + 1 < \theta + k \leq 2l + 1.
\end{aligned}
\een
Substituting these results for $\theta = 0$ in \eqref{CharIdThetaL}, one gets:
\beq
\label{FinalThetaZero}
\sum\limits_{m\in\mathbb{Z}} \frac{z^m q^{\frac{(l + 1)}{2}m(m - 1)}}{(q)_{\infty}} = \sum\limits_{i = 0}^{l} \left(z^i\,q^{i - 1 + \delta_{i, 0}} \sum\limits_{n = 0}^{\infty}\, \frac{(z\,q^{l + i})^n\,q^{\frac{l + 1}{2} n(n - 1)}}{(q)_n}\,\sum\limits_{k = 0}^{\infty} \frac{z^{-k}\,q^{\frac{l + 1}{2} k(k + 1)}}{(q)_{(l + 1)(k + \delta_{i, 0}) + i - 1}} \right),
\ee
which again can be viewed as identity in the ring $\mathbb{Z}[z^{\pm 1}][[q]]$. \eqref{FinalThetaZero} may be rewritten as
\beq
\label{AndrewsIdentityGen}
z^s: \hspace{5mm}\frac{1}{(q)_{\infty}} = \sum\limits_{\substack{n - k = s \\ n \geq 0 \\ k \geq 0}} \frac{q^{(k + s)((l + 1)k + l)}}{(q)_{k + s}(q)_{(l + 1)k + l}} + \sum\limits_{i = 0}^{l - 1} \sum\limits_{\substack{n - k + 1 = s \\ n \geq 0 \\ k \geq 0}} \frac{q^{(k + s)((l + 1)k + i)}}{(q)_{k + s - 1}(q)_{(l + 1)k + i}}
\ee
or as
\ben
z^s: \hspace{5mm}\frac{1}{(q)_{\infty}} = \sum\limits_{i = 0}^{l} \sum\limits_{\substack{n - k + 1 - \delta(i, l) = s \\ n\geq 0 \\ k\geq 0}} \frac{q^{(k + s)((l + 1)k + i)}}{(q)_{k + s - 1 + \delta(i, l)}(q)_{(l + 1)k + i}}.
\een
with fixed $s\in\mathbb{Z}$ in both forms. 
\begin{figure}[h!]
\begin{center}
\includegraphics[scale = 1]{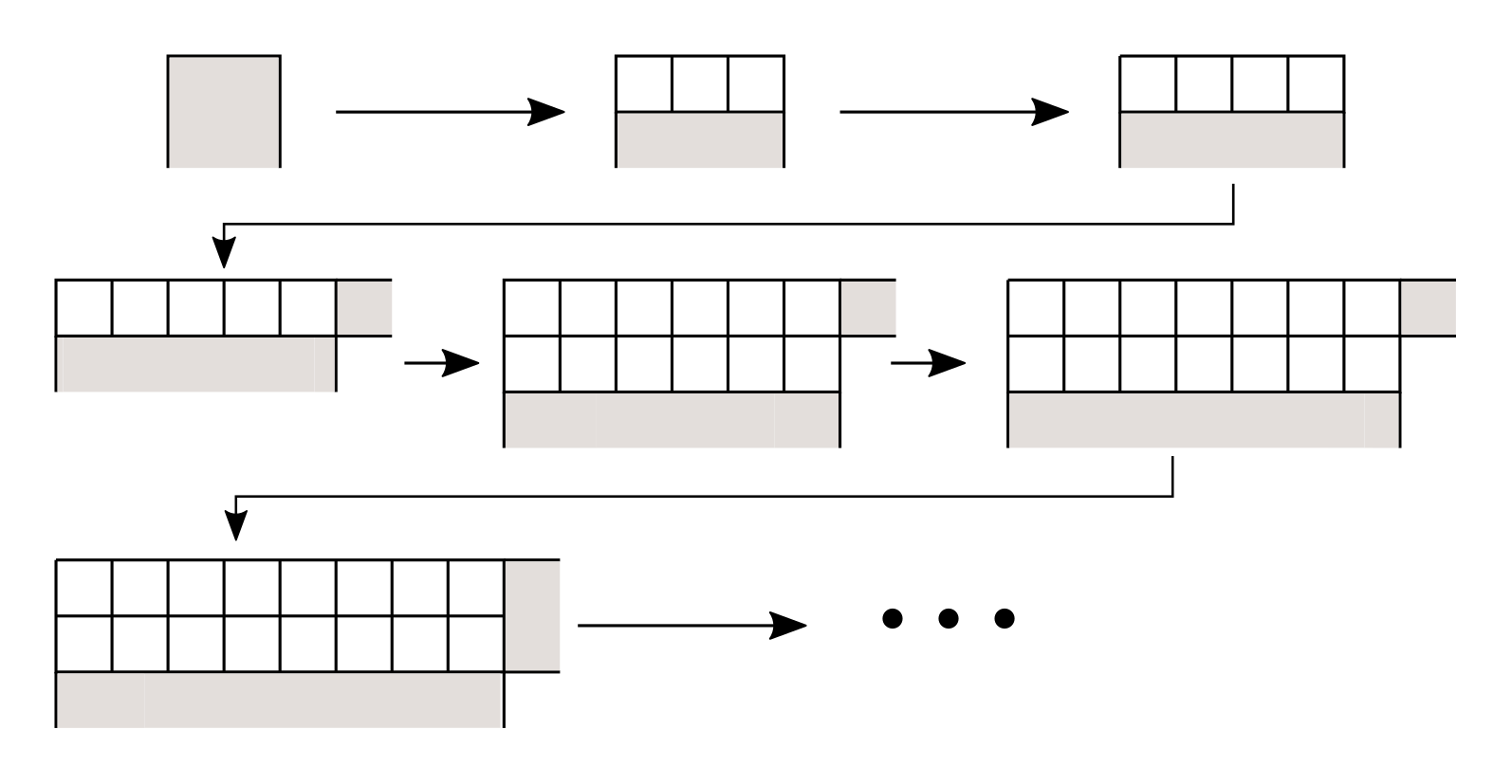}
\ben
\frac{1}{(q)_{\infty}} = \sum\limits_{k = 0}^{\infty}\frac{q^{k(3k + 2)}}{(q)_k(q)_{3k + 2}} + \sum\limits_{k = 0}^{\infty}\frac{q^{(k + 1)(3k + 3)}}{(q)_k(q)_{3k + 3}} + \sum\limits_{k = 0}^{\infty}\frac{q^{(k + 1)(3k + 4)}}{(q)_k(q)_{3k + 4}} 
\een

\caption{Identity \eqref{AndrewsIdentityGen} in the case $l = 2, \, s = 0$. \\ Durfee rectangles $k\times (3k + 2)$ \\
  Enveloping $(k + 1)\times(3k + 3)$, $(k + 1)\times(3k + 4)$}
\label{Picture 5}
\end{center}
\end{figure}
\\
General two-parametric $(n, m)$ family of this type identities could be written as
\beq
\label{Final identity l + 1}
\frac{1}{(q)_{\infty}} = \chi_{\substack{n \geq 1 \\ m\geq 1}} \sum\limits_{j = 0}^{m - 1}\frac{q^{nj}}{(q)_{n - 1}(q)_j} + \sum\limits_{k = 0}^{\infty} \frac{q^{(k + n)((l + 1)k + m)}}{(q)_{k + n}(q)_{(l + 1)k + m}} + \sum\limits_{i = 1}^{l }\sum\limits_{k = 0}^{\infty} \frac{q^{(k + n + 1)((l + 1)k + m + i)}}{(q)_{k + n}(q)_{(l + 1)k + m + i}}.
\ee
For any $n, m \in \mathbb{N}_0$, for any partition it either does not contain $n\times m$ rectangle or there is maximal $k\in\mathbb{N}_0$, such that this partition contains rectangle $(k + n)\times((l + 1)k + m)$ (such rectangle we call Durfee rectangle of this partition). Set of all partitions from the first case we denote $\tilde{P}_{l + 1, n, m}$. Generating function of $\tilde{P}_{l + 1, n, m}$ is
\ben
\chi_{\substack{n \geq 1 \\ m\geq 1}} \sum\limits_{j = 0}^{m - 1} \frac{q^{nj}}{(q)_{n - 1}(q)_j}.
\een
Denote the set of all partitions with Durfee rectangle $(k + n)\times((l + 1)k + m)$ as $P_{l + 1, k, n, m}$. Then the set of all partitions is 
\ben
P = \tilde{P}_{l + 1, n, m}\bigsqcup\left(\bigsqcup\limits_{k \in \mathbb{N}_0} P_{l + 1, k, n, m}\right).
\een
Partitions from $P_{l + 1, k, n, m}$ contain $(k + n)\times((l + 1)k + m)$ rectangle, but do not contain $(k + n)\times((l + 1)(k + 1) + m)$. Between these two rectangles there are $l$ enveloping rectangles $(k + n + 1)\times((l + 1)k + m + i)$, where $1 \leq i\leq l$. Denote the set of partitions from $P_{l + 1, k, n, m}$, containing $(k + n + 1)\times((l + 1)k + m + i)$ rectangle, but not containing $(k + n + 1)\times((l + 1)k + m + i + 1)$ rectangle $P_{l + 1, k, n, m}^{i}$. The generating function of $P_{l + 1, k, n, m}^{i}$ is
\ben
\frac{q^{(k + n + 1)((l + 1)k + i + m)}}{(q)_{k + n}(q)_{(l + 1)k + i + m}}.
\een
Thus, we get \eqref{Final identity l + 1}. Several cases with natural ordering on $P$'s are illustrated in Figures \ref{Picture 5}, \ref{Picture 6}, \ref{Picture 7} and \ref{Picture 8}. 
\\
Using formula \eqref{PInfty} one can get identities for the characters of all types of infinite Fibonacci configurations. Character identity for $\Fib_{\infty}^{(\theta,\, l)}$ is equivalent to the series of identities obtained by equating coefficients on $z^{s}$, parametrize these identities by couple $(\Fib_{\infty}^{(\theta,\, l)}, s)$. It is interesting to know if any identity of type \eqref{Final identity l + 1} could be obtained via character identity for some space of infinite Fibonacci configurations. Parametrize \eqref{Final identity l + 1} by triple $(n, m, l + 1)$. Straightforward computation gives correspondence:
\beq
\label{IdCorrespondence}
\begin{aligned}
&s \leq 0 \hspace{5mm} (\Fib_{\infty}^{(\theta,\, l)}, s) \rightarrow (0, l - s -l\,s - \theta, l + 1),
\\
&s \geq 0 \hspace{5mm} (\Fib_{\infty}^{(\theta,\, l)}, s) \rightarrow (s, l - \theta, l + 1).
\end{aligned}
\ee
One can verify this correspondence in the particular case of $l = 1$ with \eqref{11pos}, \eqref{11neg}, \eqref{01pos}, \eqref{01neg}. Identities \eqref{Final identity l + 1} with fixed $l$ may be viewed as the two-parametric deformation of identity $(0, 0, l + 1)$ derived in \cite{Andrews_Part1}. With fixed $l$ identities $(n, m, l + 1)$ may be identified with points $\mathbb{N}_0\times\mathbb{N}_0 \subset \mathbb{R}^2$ of $(n, m)$ plane. Then identities $(\Fib_{\infty}^{(\theta,\, l)}, s)$ lie exactly on lines $n = 0$, $m = i$, $0 \leq i \leq l$, see Figure \ref{Diagram}.
\begin{remark}
\label{RemarkLine}
Identities \eqref{Final identity l + 1} for parameters $(n_1, m_1)$ and $(n_2, m_2)$ lying on some line of type $m = (l + 1)\,n + m'$ are equivalent. Indeed, let $n_2 = n_1 + 1$, $m _2 = m_1 + (l + 1)$. Identity for $(n_1, m_1)$
\ben
\frac{1}{(q)_{\infty}} = \chi_{\substack{n_1 \geq 1 \\ m_1\geq 1}} \sum\limits_{j = 0}^{m_1 - 1}\frac{q^{n_1\,j}}{(q)_{n_1 - 1}(q)_j} + \sum\limits_{k = 0}^{\infty} \frac{q^{(k + n_1)((l + 1)k + m_1)}}{(q)_{k + n_1}(q)_{(l + 1)k + m_1}} + \sum\limits_{i = 1}^{l}\sum\limits_{k = 0}^{\infty} \frac{q^{(k + n_1 + 1)((l + 1)k + m_1 + i)}}{(q)_{k + n_1}(q)_{(l + 1)k + m_1 + i}}
\een  
may be rewritten as
\ben
\begin{aligned}
&\frac{1}{(q)_{\infty}} = \chi_{\substack{n_1 \geq 1 \\ m_1\geq 1}} \sum\limits_{j = 0}^{m_1 - 1}\frac{q^{n_1\,j}}{(q)_{n_1 - 1}(q)_j} + \frac{q^{n_1\,m_1}}{(q)_{n_1}(q)_{m_1}} + \sum\limits_{i = 1}^{l} \frac{q^{((n_1 + 1)(m_1 + i))}}{(q)_{n_1}(q)_{m_1 + i}}
\\
&+ \sum\limits_{k = 0}^{\infty} \frac{q^{(k + n_2)((l + 1)k + m_2)}}{(q)_{k + n_2}(q)_{(l + 1)k + m_2}} + \sum\limits_{i = 1}^{l}\sum\limits_{k = 0}^{\infty} \frac{q^{(k + n_2 + 1)((l + 1)k + m_2 + i)}}{(q)_{k + n_2}(q)_{(l + 1)k + m_2 + i}}. 
\end{aligned}
\een 
Then we need to check 
\ben
\chi_{\substack{n_1 \geq 1 \\ m_1\geq 1}} \sum\limits_{j = 0}^{m_1 - 1}\frac{q^{n_1\,j}}{(q)_{n_1 - 1}(q)_j} + \frac{q^{n_1\,m_1}}{(q)_{n_1}(q)_{m_1}} + \sum\limits_{i = 1}^{l} \frac{q^{((n_1 + 1)(m_1 + i))}}{(q)_{n_1}(q)_{m_1 + i}}
\een
equals
\ben
\sum\limits_{j = 0}^{m_2 - 1}\frac{q^{n_2\,j}}{(q)_{n_2 - 1}(q)_j} = \sum\limits_{j = 0}^{m_1 + l}\frac{q^{(n_1 + 1)\,j}}{(q)_{n_1}(q)_j},
\een
or equivalently
\beq
\label{supportIdentity}
\chi_{\substack{n_1 \geq 1 \\ m_1\geq 1}} \sum\limits_{j = 0}^{m_1 - 1}\frac{q^{n_1\,j}}{(q)_{n_1 - 1}(q)_j} + \frac{q^{n_1\,m_1}}{(q)_{n_1}(q)_{m_1}} = \sum\limits_{j = 0}^{m_1}\frac{q^{(n_1 + 1)\,j}}{(q)_{n_1}(q)_j}, 
\ee
which is true as soon as both parts are characteristic functions of partitions not containing $(n_1 + 1)\times(m_1 + 1)$ rectangle. 
Thus, all identities \eqref{Final identity l + 1} can be obtained via character identity for $\Fib_{\infty}^{(\theta,\, l)}$ with respect to \eqref{IdCorrespondence}.
\end{remark}

\section{Lattice vertex algebras}
\label{Section 3}
Heisenberg Lie algebra $\mathfrak{H}$ is $\mathfrak{H} = \langle C, a_n \:|\: n\in\mathbb{Z}  \rangle_{\mathbb{C}}$ with bracket
\ben
[a_n, a_m] = \delta_{n + m, 0}\:C, \hspace{5mm} [C, \:\cdot\:] = 0. 
\een
Let $F_{\mu}$ be the Fock module over Heisenberg algebra with the highest vector $|\mu\rangle$:
\ben
c |\mu\rangle = |\mu\rangle,\hspace{5mm} h_0 |\mu\rangle = \mu|\mu\rangle.  
\een
Vertex operator superalgebra structure on $V_{\sqrt{N}\,\mathbb{Z}}\, = \bigoplus\limits_{\mu \in \sqrt{N}\,\mathbb{Z}}\, F_{\mu}$ is defined by follows:
\begin{enumerate}
    \item $(\mathbb{Z}/2\mathbb{Z})$\textbf{-gradation}:
    \ben
    V^{\Bar{0}} : = \bigoplus\limits_{\substack{\lambda\in\sqrt{N}\mathbb{Z} \\ \lambda^2 \equiv 0 \mod 2}}\,F_{\lambda},\hspace{1cm} V^{\Bar{1}} : = \bigoplus\limits_{\substack{\lambda\in\sqrt{N}\mathbb{Z} \\ \lambda^2 \equiv 1 \mod 2}}\,F_{\lambda}.
    \een
    \item (\textbf{Vacuum vector}) \, $|0\rangle$.
    \item (\textbf{Translation operator}) \, $T = \sum\limits_{n \geq 0}\, a_{- n - 1}\,a_n$. 
    \item (\textbf{Vertex operators}) Set
    \ben
\begin{aligned}
&Y(|0\rangle, z) = 0, \hspace{5mm} Y(a_{-1}\,|0\rangle) := a(z) = \sum\limits_{n\in\mathbb{Z}}\,a_n\,z^{- n - 1}.
\\
&Y(|\mu\rangle, z) = e^{\mu\, q}\,z^{\mu\,a_0}\,\exp\left( -\mu\sum\limits_{n < 0}\, \frac{a_{n}}{n}\,z^{-n}\right)\,\exp\left( -\mu\sum\limits_{n > 0}\, \frac{a_{n}}{n}\,z^{-n}\right),
\end{aligned}
    \een
\end{enumerate}
with $[a_0, q] = 1$. Due to the Strong reconstruction theorem these data are enough to define a vertex superalgebra structure on $V_{\sqrt{N}\,\mathbb{Z}}$ (see \cite{Dong_Lepowsky, Frenkel2004, Iohara} for details).
\\\\
OPE of two vertex operators is 
\ben
V_{\eta}(z)\,V_{\mu}(\omega) = (z - w)^{\eta\mu}\, :V_{\eta}(z)\,V_{\mu}(\omega):,
\een
where
\ben
:V_{\eta}(z)\,V_{\mu}(\omega): = e^{(\eta + \mu)\,q}\,z^{\eta\,a_0}\,w^{\mu\,a_0}\,\exp\left( -\sum\limits_{n < 0}\, \frac{a_{n}}{n}\,(\eta\,z^{-n} + \mu\,w^{-n})\right)\,\exp\left( -\sum\limits_{n > 0}\, \frac{a_{n}}{n}\,(\eta\,z^{-n} + \mu\,w^{-n})\right).
\een
There is a family of conformal vectors $\omega_{\lambda} = \frac{1}{2}\,a_{-1}^2 + \lambda\,a_{-2}$, $\lambda \in\mathbb{C}$ in $V_{\sqrt{N}\,\mathbb{Z}}$. In this work vertex algebra $V_{\sqrt{N}\mathbb{Z}}$ is considered as conformal with conformal vector $\frac{1}{2}\,a_{-1}^2 + \frac{\sqrt{N}}{2}\,a_{-2}$. Denote
\ben
V_{(i), \sqrt{N}} = \bigoplus\limits_{m \in\mathbb{Z}} F_{\frac{i}{\sqrt{N}} + m\sqrt{N}} \hspace{5mm} i = 0, 1, \ldots, N - 1
\een
the set of irreducible representations of $V_{\sqrt{N}\mathbb{Z}}$. These modules are bigraded by the operator $L_0$ and charge $m$. Define the character as
\ben
\ch V_{(i), \sqrt{N}} = \Tr\left(z^{\frac{h_0}{\sqrt{N}} - \frac{i}{N}}\,q^{L_0}|_{V_{(i), \sqrt{N}}}\right) = q^{\frac{i^2}{2N} - \frac{i}{2}}\,\sum\limits_{m\in\mathbb{Z}}\,\frac{(z\,q^i)^m\,q^{N\,\frac{m(m - 1)}{2}}}{(q)_{\infty}}.
\een
Denote $\Fib_{\infty} = \bigcup\limits_{0 \leq \theta \leq l \in\mathbb{N}_0} \,\Fib_{\infty}^{(\theta,\, l)}$.
Construct by the induction map $\tau\colon\Fib_{\infty} \longrightarrow \mathbb{Z}^{\mathbb{N}}$ as
\ben
\tau(a)_i = -\max \left(\{j\in\mathbb{Z}\hspace{2mm} |\hspace{2mm} a_j\neq 0\}\setminus\{\tau(a)_{i - 1}, \ldots, \tau(a)_{1}\}\right)
\een
This map assigns to any Fibonacci configuration $a(n)$ the increasing sequence of positions with ``1" of the configuration $a(-n)$. Then \hyperref[TheoremBasis]{Theorem} \ref{TheoremBasis} can be reformulated~---~the set
\ben
\{\mathbf{e}_a = e_{\tau(a)_1}\,e_{\tau(a)_2}\,e_{\tau(a)_3}\dots \hspace{2mm}|\hspace{2mm} a\in\Fib^{(i, N)}_{\infty}\},
\een
is the basis in $V_{(i), \sqrt{N}}$. Notice that the characters of $V_{(i), \sqrt{N}}$ and $\Fib^{(i, N)}_{\infty}$ differs by multiple $q^{\frac{i^2}{2N} - \frac{i}{2}}$, i.e.\ character identities derived in the last section are character identities for $V_{(i), \sqrt{N}}$.  
\\\\
Semi-infinite configurations $\Fib_{\infty}^{(\theta,\, l)}$ appear in context of rank one lattice vertex operator (super)algebras (and their modules) with lattice parameter $\langle\lambda|\lambda\rangle = l + 1\in \mathbb{N}$ and correspond to the basis of semi-infinite monomials, constructed by the Fourier modes of vertex operator $V_{\lambda}(z)$. Character identities turn out to be Durfee rectangle identities with the ratio of height to base $1 : l + 1$. Nevertheless, G. Andrews in \cite{Andrews_Part1} describes Durfee rectangle identities with the ratio of height to base being an	 arbitrary positive rational number. This observation leads to
\begin{conjecture}
Suppose $p, q \in\mathbb{N}, (p, q) = 1$. Let the space
\ben
V_{\sqrt{\frac{p}{q}}\,\mathbb{Z}} = \bigoplus\limits_{m\in\mathbb{Z}}\, F_{m\sqrt{\frac{p}{q}}}
\een
be endowed with action of bosonic vertex operator $V_{\sqrt{\frac{p}{q}}}(z)$. Then $V_{\sqrt{\frac{p}{q}}\,\mathbb{Z}}$ admits Feigin-Stoyanovsky basis in terms of $V_{\sqrt{\frac{p}{q}}}(z)$ modes, whose character identity is equivalent to series of $p/q$ Durfee rectangle identity (\cite{Andrews_Part1}) with shifts.
\end{conjecture}

\clearpage
\bibliography{bibliography}{}

\begin{thebibliography}{CLM08b}

\bibitem[And71]{Andrews_Part1}
George~E. Andrews.
\newblock {Generalizations of the Durfee Square}.
\newblock {\em Journal of the London Mathematical Society}, s2-3(3):563--570, 1971.

\bibitem[And84]{AndrewsBook}
George~E. Andrews.
\newblock {\em {The Theory of Partitions}}.
\newblock Encyclopedia of Mathematics and its Applications. Cambridge University Press, 1984.

\bibitem[Bar11]{Baranovi2011}
Ivana Baranovi{\'{c}}.
\newblock {Combinatorial Bases of Feigin{\textendash}Stoyanovsky's Type Subspaces of Level 2 Standard Modules for}.
\newblock {\em Communications in Algebra}, 39(3):1007--1051, March 2011.

\bibitem[Bre89]{Bressoud1989}
David~M. Bressoud.
\newblock {Lattice paths and the Rogers-Ramanujan identities}.
\newblock In {\em Lecture Notes in Mathematics}, pages 140--172. Springer Berlin Heidelberg, 1989.

\bibitem[Cal07]{Calinescu2007}
Corina Calinescu.
\newblock {Principal subspaces of higher-level standard $\widehat{sl(3)}$-modules}.
\newblock {\em Journal of Pure and Applied Algebra}, 210(2):559--575, August 2007.

\bibitem[CLM08a]{CALINESCU2008}
C.~Calinescu, J.~Lepowsky, and A.~Milas.
\newblock {Vertex-algebraic structure of the principal subspaces of certain $A_1^{(1)}$-modules. {I}: Level one case.}
\newblock {\em International Journal of Mathematics}, 19(01):71--92, January 2008.

\bibitem[CLM08b]{Calinescuu2008}
C.~Calinescu, J.~Lepowsky, and A.~Milas.
\newblock {Vertex-algebraic structure of the principal subspaces of certain $A_1^{(1)}$-modules. {II}: Higher-level case}.
\newblock {\em Journal of Pure and Applied Algebra}, 212(8):1928--1950, August 2008.

\bibitem[CLM10]{Calinescu2010}
C.~Calinescu, J.~Lepowsky, and A.~Milas.
\newblock {Vertex-algebraic structure of the principal subspaces of level one modules for the untwisted affine Lie algebras of types $A, D, E$}.
\newblock {\em Journal of Algebra}, 323(1):167--192, January 2010.

\bibitem[DL12]{Dong_Lepowsky}
Chongying Dong and James Lepowsky.
\newblock {\em {Generalized Vertex Algebras and Relative Vertex Operators}}.
\newblock Birkhäuser Boston, MA, 2012.

\bibitem[FBZ04]{Frenkel2004}
Edward Frenkel and David Ben-Zvi.
\newblock {\em Vertex Algebras and Algebraic Curves}.
\newblock American Mathematical Society, August 2004.

\bibitem[FF06]{FF_Principle}
B.~Feigin and E.~Feigin.
\newblock {Principal subspace for the bosonic vertex operator $\phi_{\sqrt{2m}}(z)$ and Jack polynomials}.
\newblock {\em Advances in Mathematics}, 206(2):307--328, November 2006.

\bibitem[IK11]{Iohara}
Kenji Iohara and Yoshiyuki Koga.
\newblock {\em {Representation Theory of the Virasoro Algebra}}, volume 308.
\newblock Springer Monographs in Mathematics, 01 2011.

\bibitem[Jac29]{Jacobi_Fundamenta}
C.~G. Jacobi.
\newblock {\em {Fundamenta nova theoriae functionum ellipticarum}}.
\newblock Sumtibus fratrum, 1829.

\bibitem[MP12]{Penn_Lattice}
Antun Milas and Michael Penn.
\newblock {Lattice vertex algebras and combinatorial bases: General case and W-algebras}.
\newblock {\em The New York Journal of Mathematics [electronic only]}, 18:621 -- 650, 08 2012.

\bibitem[Rog94]{Rogers1894}
L.~J. Rogers.
\newblock Third memoir on the expansion of certain infinite products.
\newblock {\em Proceedings of the London Mathematical Society}, s1-26(1):15--32, November 1894.

\bibitem[RR19]{RogRam}
S.~Ramanujan and L.~J. Rogers.
\newblock {Proof of certain identities in combinatory analysis}.
\newblock {\em Cambr. Phil. Soc. Proc.}, 19:211--216, 1919.

\bibitem[Sch73]{Schur}
I.~Schur.
\newblock Ein beitrag zur additiven zahlentheorie und zur theorie der kettenbr¨uch.
\newblock {\em S.-B. Preuss. Akad. Wiss. Phys.-Math. Kl.}, 2, 01 1973.

\bibitem[SF94]{FS}
A.~V. Stoyanovsky and B.~L. Feigin.
\newblock {Functional models for representations of current algebras and semi-infinite Schubert cells}.
\newblock {\em Functional Analysis and Its Applications}, 28(1):55--72, 1994.

\bibitem[Tru10]{Trupevi2010}
Goran Trup{\v{c}}evi{\'{c}}.
\newblock {Combinatorial Bases of Feigin{\textendash}Stoyanovsky's Type Subspaces of Level 1 Standard Modules for}.
\newblock {\em Communications in Algebra}, 38(10):3913--3940, September 2010.

\bibitem[Tru11]{zbMATH05919368}
G.~Trup{\v{c}}evi{\'c}.
\newblock {Characters of {Feigin}-{Stoyanovsky}'s type subspaces of level one modules for affine {Lie} algebras of types {{\(A_{\ell}^{(1)}\)}} and {{\(D_4^{(1)}\)}}}.
\newblock {\em Glas. Mat., III. Ser.}, 46(1):49--70, 2011.

\bibitem[War23]{Warnaar2023}
S.~O. Warnaar.
\newblock {The $A_2$ Andrews{\textendash}Gordon identities and cylindric partitions}.
\newblock {\em Transactions of the American Mathematical Society, Series B}, 10(22):715--765, June 2023.

\end{thebibliography}

\clearpage
\appendix
\section{Images}
\begin{figure}[h!]
\begin{center}
\includegraphics[scale = 1.1]{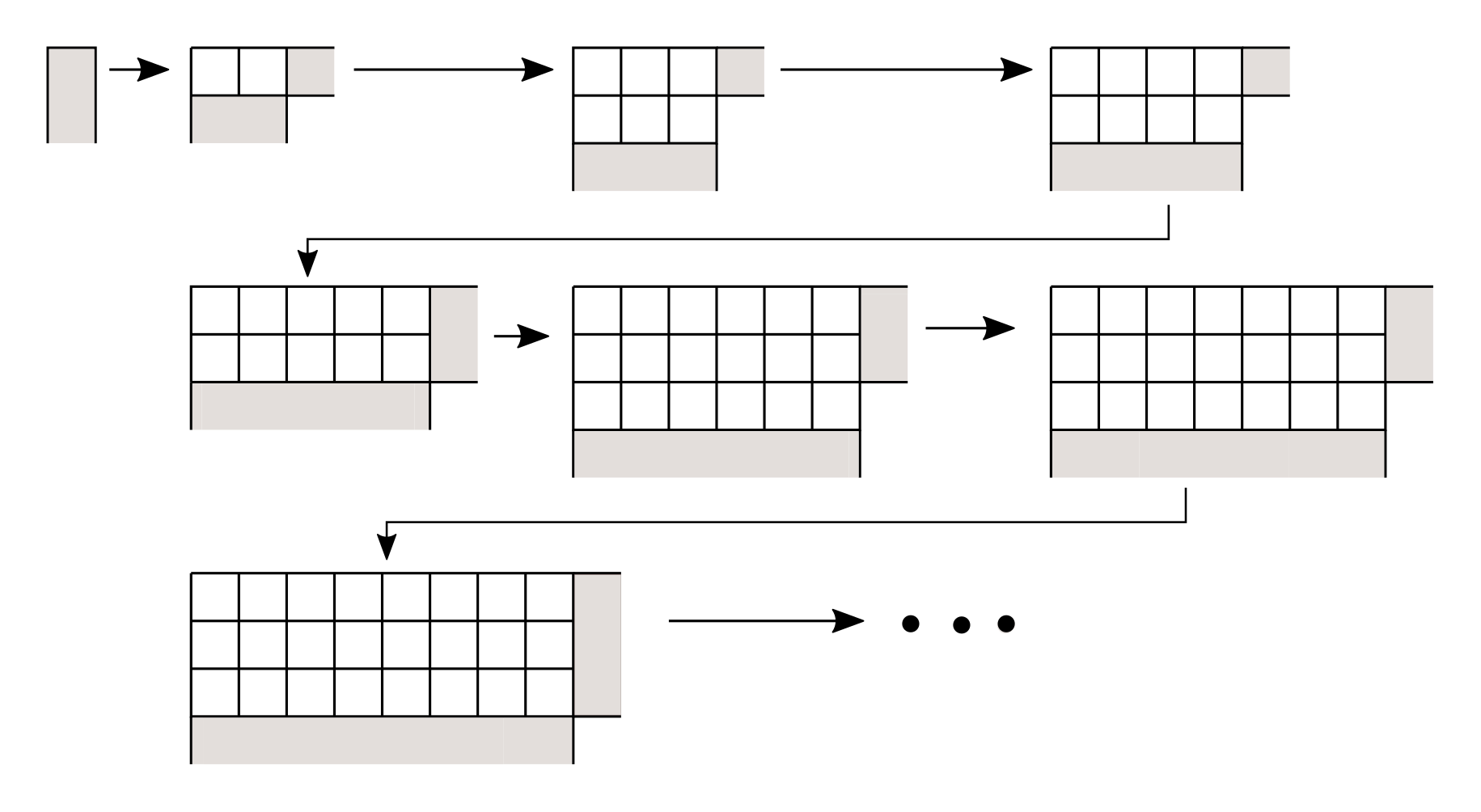}
\ben
\frac{1}{(q)_{\infty}} = \sum\limits_{k = 0}^{\infty}\frac{q^{(k + 1)(3k + 2)}}{(q)_{k + 1}(q)_{3k + 2}} + \sum\limits_{k = 0}^{\infty}\frac{q^{(k + 1)3k}}{(q)_k(q)_{3k}} + \sum\limits_{k = 0}^{\infty}\frac{q^{(k + 1)(3k + 1)}}{(q)_k(q)_{3k + 1}} 
\een

\caption{Identity \eqref{AndrewsIdentityGen} in case $l = 2, \, s = 1$. \\ Durfee rectangles: $(k + 1)\times (3k + 2)$ \\
  Enveloping rectangles: $(k + 1)\times3k$, $(k + 1)\times(3k + 1)$}
\label{Picture 6}
\end{center}
\end{figure}

\begin{figure}[h!]
\begin{center}
\includegraphics[scale = 1.1]{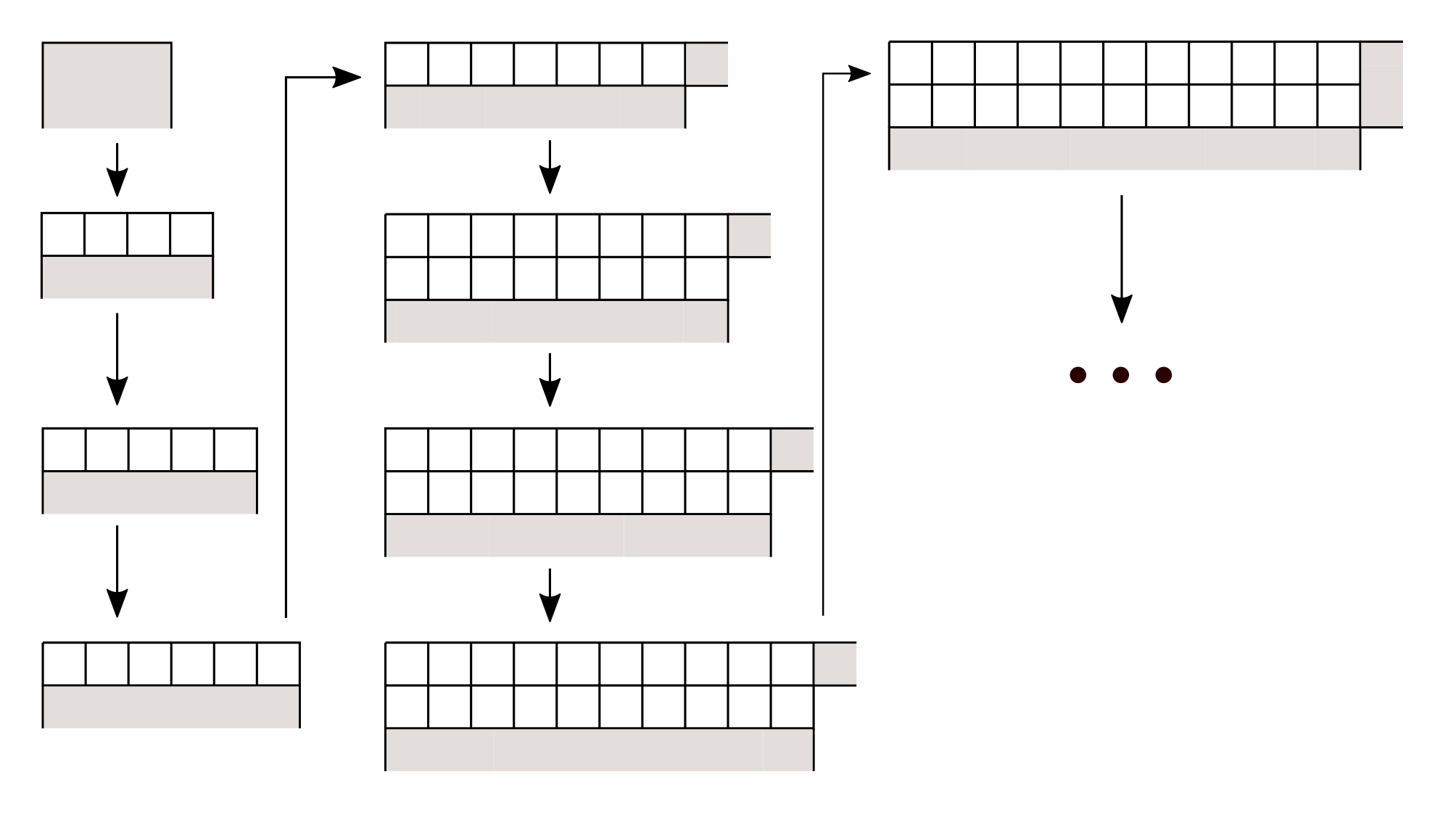}
\ben
\frac{1}{(q)_{\infty}} = \sum\limits_{k = 0}^{\infty} \frac{q^{k(4k + 3)}}{(q)_k (q)_{4k + 3}} + \sum\limits_{k = 0}^{\infty} \frac{q^{(k + 1)(4k + 4)}}{(q)_k (q)_{4k + 4}} + \sum\limits_{k = 0}^{\infty} \frac{q^{(k + 1)(4k + 5)}}{(q)_k (q)_{4k + 5}} + \sum\limits_{k = 0}^{\infty} \frac{q^{(k + 1)(4k + 6)}}{(q)_k (q)_{4k + 6}} 
\een

\caption{Identity \eqref{AndrewsIdentityGen} in case $l = 3, \, s = 0$. \\ Durfee rectangles: $k\times (4k + 3)$ \\
  Enveloping rectangles: $(k + 1)\times(4k + 4)$, $(k + 1)\times(4k + 5)$, $(k + 1)\times(4k + 6)$}
\label{Picture 7}
\end{center}
\end{figure}

\begin{figure}[h!]
\begin{center}
\includegraphics[scale = 1.1]{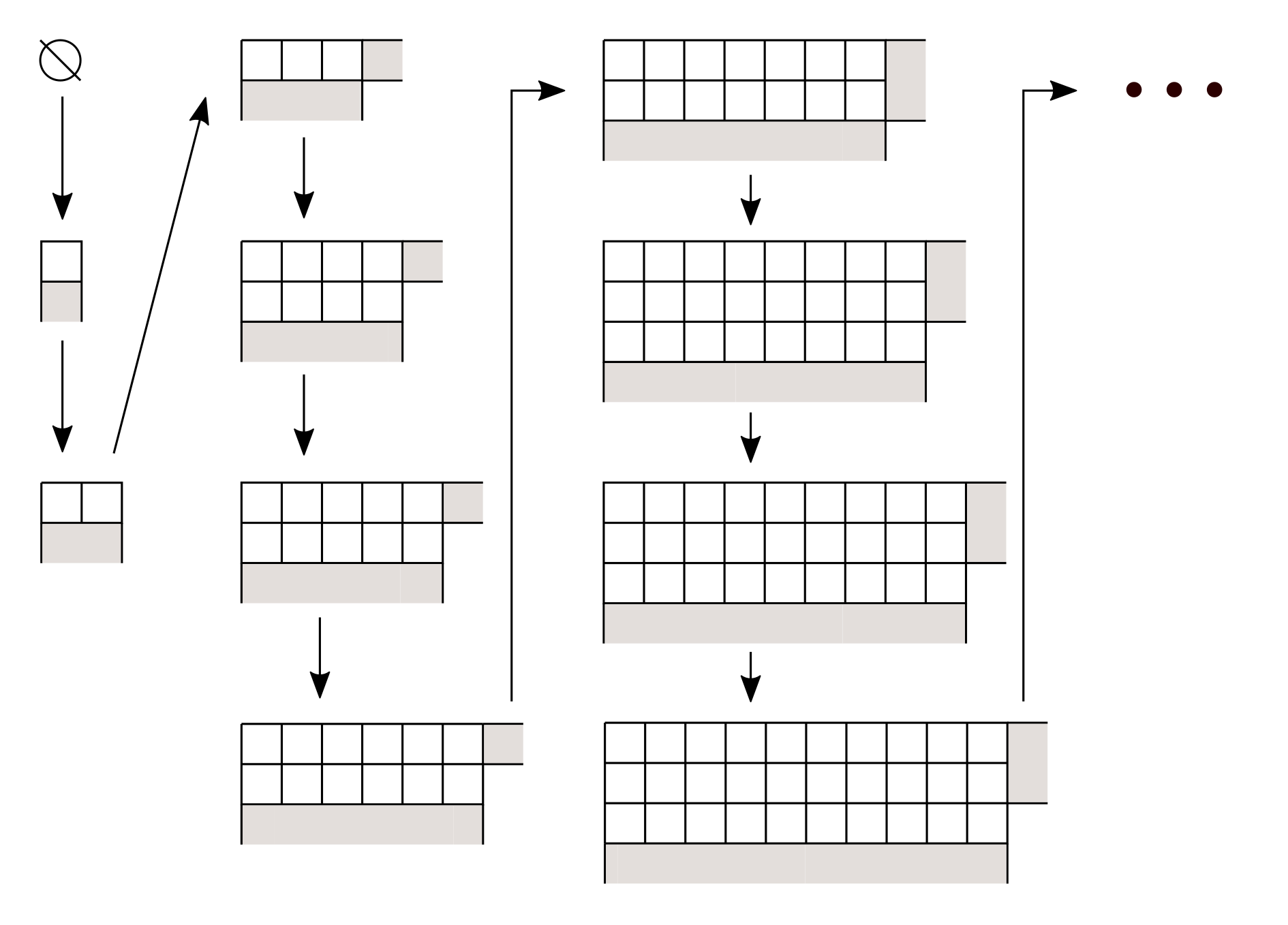}
\ben
\frac{1}{(q)_{\infty}} = \sum\limits_{k = 0}^{\infty} \frac{q^{(k + 1)(4k + 3)}}{(q)_{k + 1} (q)_{4k + 3}} + \sum\limits_{k = 0}^{\infty} \frac{q^{(k + 1)4k}}{(q)_k (q)_{4k}} + \sum\limits_{k = 0}^{\infty} \frac{q^{(k + 1)(4k + 1)}}{(q)_k (q)_{4k + 1}} + \sum\limits_{k = 0}^{\infty} \frac{q^{(k + 1)(4k + 2)}}{(q)_k (q)_{4k + 2}} 
\een

\caption{Identity \eqref{AndrewsIdentityGen} in case $l = 3, \, s = 1$. \\ Durfee rectangles: $(k + 1)\times (4k + 3)$ \\
  Enveloping rectangles: $(k + 1)\times 4k $, $(k + 1)\times(4k + 1)$, $(k + 1)\times(4k + 2)$}
\label{Picture 8}
\end{center}
\end{figure}

\end{document}